\documentstyle[12pt,amssym]{article}
\newtheorem{theo}{Theorem}
\newtheorem{prop}[theo]{Proposition}
\newtheorem{lemm}[theo]{Lemma}
\newtheorem{coro}[theo]{Corollary}

\newcommand{\myF}{{\cal F}}

\newcommand{\myB}{{\cal B}}

\newcommand{\myg}{{\frak g}}

\newcommand{\myk}{{\frak k}}

\newcommand{\myp}{{\frak p}}

\newcommand{\bR}{{\Bbb R}}

\newcommand{\Ad}{{\rm Ad}}

\newcommand{\trace}{{\rm Trace}}

\newcommand{\id}{{\rm id}}

\newcommand{\levy}{L\'{e}vy }
\newcommand{\lims}{\overline{\lim}}

\begin{document}

\begin{center}
{\large \bf Inhomogeneous \levy processes in Lie groups and

homogeneous spaces}

Ming Liao \ (Auburn University)
\end{center}

\begin{quote}
{\bf Summary} \ We obtain a representation of an inhomogeneous \levy
process in a Lie group or a homogeneous space in terms of a drift, a matrix function
and a measure function. Because the stochastic continuity is not assumed,
our result generalizes the well known L\'{e}vy-It\^{o} representation
for stochastic continuous processes with independent increments
in $\bR^d$ and its extension to Lie groups.



\end{quote}

\section{Introduction} \label{intro}

Let $x_t$, $t\in\bR_+=[0,\,\infty)$, be a process in $\bR^d$ with rcll paths
(right continuous paths with left limits). It is said to
have independent increments if for $s<t$, $x_t-x_s$ is independent
of $\myF_s^x$ (the $\sigma$-algebra generated by $x_u$, $0\leq u
\leq s$). The process is called a \levy process if it also has stationary increments,
that is, if $x_t-x_s$ has the same distribution as $x_{t-s}-x_0$.
It is well known that the class of \levy processes
in $\bR^d$ coincides with the class of rcll Markov processes with
translation-invariant transition functions $P_t$. The celebrated
\levy-Khinchin formula provides a useful representation of a \levy
process in $\bR^d$ by a triple of a drift vector, a covariance matrix and a \levy measure.

More generally, a rcll process in $\bR^d$ with independent but possibly non-stationary increments may be called an inhomogeneous \levy process.
It is easy to show that such processes, also called additive processes in literature, coincide with rcll inhomogeneous Markov processes
with translation invariant two-parameter transition functions $P_{s,t}$ (see \cite{sato}).

A \levy process $x_t$ is always stochastically continuous, that is, $x_t=x_{t-}$ a.s. (almost surely)
for any fixed $t>0$. An inhomogeneous \levy process $x_t$ may not be stochastically continuous,
but if it is, then the well known L\'{e}vy-It\^{o} representation holds (\cite[chapter~15]{kallenberg}):
\begin{equation}
x_t = x_0 + b_t + B_t + \int_{[0,\,t]}\int_{\vert x\vert\leq
1}x\tilde{N}(ds,dx) + \int_{[0,\,t]}\int_{\vert x\vert>1}x N(ds,dx),
\label{levyitoRn}
\end{equation}
where $b_t$ is a continuous path in $\bR^d$ with $b_0=0$, called a drift, $B_t=(B_t^1,\ldots,B_t^d)$ is a $d$-dim continuous Gaussian process
of zero mean and independent increments, $N(dt,dx)$ is an independent Poisson random measure on $\bR_+\times\bR^d$ with
intensity $\eta=E(N)$ being a \levy measure function (to be defined later), and $\tilde{N}=N-\eta$ is the compensated form of $N$.

The distribution of the Gaussian process $B_t$ is determined by its covariance matrix
function $A_{ij}(t)=E(B_t^iB_t^j)$ and the distribution of the Poisson random measure $N$ is determined
by its intensity measure $\eta$. Thus, the distribution of a stochastically continuous
inhomogeneous \levy process $x_t$ in $\bR^d$ is completely determined by the time-dependent triple $(b,A,\eta)$.

In general, an inhomogeneous \levy process $x_t$ may not be a semimartingale (see \cite[II]{js}). By It\^{o}'s formula, it can be shown that the
process $z_t=x_t-b_t$ is a semi-martingale and for any smooth function $f(x)$ on $\bR^d$ of compact support,
\begin{equation} f(z_t) -
\frac{1}{2}\sum_{i,j}\int_0^tf_{ij}(z_s)dA_{ij}(s) -
\int_0^t\int_{\bR^d}[f(z_s+x)-f(z_s)-\sum_ix_if_i(z_s)1_{[\vert
x\vert\leq 1]}]\eta(ds,dx) \label{martRn}
\end{equation}
is a martingale, where $f_i$ and $f_{ij}$ denote respectively the
first and second order partial derivatives of $f$. This in fact
provides a complete characterization for the distribution
of a stochastically continuous inhomogeneous \levy process $x_t$ in $\bR^d$.

This martingale representation is extended to stochastically continuous inhomogeneous \levy process in a general Lie group $G$ in \cite{feinsilver},
generalizing an earlier result in \cite{stroockvaradhan} for continuous processes. A different form of martingale representation
in terms of the abstract Fourier analysis is obtained in \cite{heyermartch}, where the processes considered are also stochastically continuous.

A \levy process $x_t$ in a Lie group $G$ is defined as a rcll process with independent and stationary (multiplicative) increments, that is,
for $s<t$, $x_s^{-1}x_t$ is independent of $\myF_s^x$ and has the same distribution as $x_0^{-1}x_{t-s}$. Such a process may also
be characterized as a rcll Markov process in $G$ whose transition function $P_t$ is invariant under left translations on $G$.
The classical triple representation of \levy processes was extended to Lie groups in \cite{hunt} in the form of a generator formula. A
functional form of L\'{e}vy-It\^{o} representation for \levy processes in Lie groups was obtained in \cite{ak}.

An inhomogeneous \levy process $x_t$ in a Lie group $G$ is defined to be a rcll process that has independent but not necessarily stationary
increments, which may also be characterized as an inhomogeneous rcll Markov process with a left invariant transition function $P_{s,t}$ (see
Proposition~\ref{prinhlevyG}). As mentioned earlier, the stochastically continuous
inhomogeneous \levy processes in $G$ may be represented by a martingale determined by a triple $(b,A,\eta)$.

The notion of \levy processes as invariant Markov processes, including inhomogeneous ones, may be extended to more general
homogeneous spaces, such as a sphere. A homogeneous space $G/K$ may be regarded as a manifold $X$ under the transitive action of a Lie group $G$
with $K$ being the subgroup that fixes a point in $X$. As on a Lie group $G$, a Markov process in $X$ with a $G$-invariant transition function $P_t$
will be called a \levy process, and an inhomogeneous Markov process with a $G$-invariant transition function $P_{s,t}$ will be called
an inhomogeneous \levy process. Although there is no natural product on $X=G/K$, the increments of a process in $X$ may
be properly defined, and a \levy process in $X$ may be characterized by independent and stationary increments, with inhomogeneous ones just by
independent increments, same as on $G$. See section \S\ref{sec6} for more details.

The purpose of this paper is to study the representation of inhomogeneous \levy processes,
not necessarily stochastically continuous, in a Lie group $G$ or a homogeneous space $G/K$.
We will show that such a process is represented by a triple $(b,A,\eta)$ with possibly discontinuous $b_t$
and $\eta(t,\cdot)=\eta([0,\,t]\times\cdot)$. This is in contrast to the stochastically continuous case when the triple is continuous in $t$.
A non-stochastically-continuous process may have a fixed jump, that is, a time $t>0$ such that $P(x_t\neq x_{t-})>0$.
It turns out to be quite non-trivial to handle fixed jumps which may form a countable dense subset of $\bR_+$.

Our result applied to $\bR^n$ leads to a martingale representation of inhomogeneous \levy processes in $\bR^n$, for which
the martingale (\ref{martRn}) contains an extra term:
\[-\sum_{u\leq t}\int_{\bR^n}[f(z_{u-}+x-h_u)-f(z_{u-})]\nu_u(dx),\]
where $\nu_u$ is the distribution of the fixed jump at time $u$ with mean $h_u$.
This complements the Fourier transform representation on $\bR^n$ obtained in \cite[II.5]{js}.

We will obtain the representation of inhomogeneous \levy processes not only on a Lie group $G$ but also on a homogeneous space $G/K$. For this purpose,
we will formulate a product structure and develop certain invariance properties on $G/K$ so that the formulas obtained on $G$, and their proofs, may
be carried over to $G/K$. We will also show that an inhomogeneous \levy process in $G/K$ is the natural projection of an inhomogeneous \levy
process in $G$. On an irreducible $G/K$, such as a sphere, the representation takes a very simple form: there is no drift and
the covariance matrix function $A(t)=a(t)I$ for some function $a(t)$. This is even simpler than the representation on $\bR^n$.

Our interest in non-stochastically-continuous inhomogeneous \levy
processes lies in the following application. Let $x_t$ be a Markov
process in a manifold $X$ invariant under the action of a Lie group $G$.
It is shown in \cite{Liao2009} that $x_t$ may be decomposed
into a radial part $y_t$, transversal to $G$-orbits, and
an angular part $z_t$, in a standard $G$-orbit $Z$.
The $y_t$ can be an arbitrary Markov process in a transversal subspace,
whereas given $y_t$, the conditioned $z_t$ is
an inhomogeneous \levy process in the homogeneous space $Z=G/K$. For example, a Markov process
in $\bR^n$ invariant under the group $O(n)$ of orthogonal
transformations is decomposed into a radial Markov process in a
half line and an angular process in the unit sphere.

In \cite{Liao2009}, the representation of stochastically continuous inhomogeneous \levy
processes in $G/K$ is used to obtain a skew product decomposition of a $G$-invariant continuous Markov
process $x_t$ in which the angular part $z_t$ is a time changed Brownian motion in $G/K$, generalizing
the well known skew product of Brownian motion in $\bR^n$.

When the $G$-invariant Markov process $x_t$ is discontinuous, its conditioned angular part $z_t$ is typically
not stochastically continuous. For example,
a discontinuous $O(n)$-invariant \levy process in $\bR^n$ is stochastically continuous,
but its conditioned angular part is not. The present result will provide a useful tool
in this situation.

We note that an important related problem, the weak convergence of convolution
products of probability measures to a two-parameter convolution semigroup, which is the distribution
of an inhomogeneous \levy process, is not pursued in this paper. The stochastically continuous case is studied
in \cite{siebert} and \cite{pap}.

Our paper is organized as follows. The main results on Lie groups are stated in the next section with proofs given in the four sections to follow.
In \S\ref{secmeasfun},
we establish the martingale representation, associated to a triple $(b,A,\eta)$, for a given inhomogeneous \levy process,
under two technical assumptions (A) and (B). These assumptions are verified in \S\ref{secAB}. We then prove the uniqueness of the triple
for a given process in \S\ref{secuniquetriple}, and the uniqueness and the existence of the process for a given triple in \S\ref{secuniqueexistprocess}.
The results on homogeneous spaces are presented in \S\ref{sec6}.
We follow the basic ideas in \cite{feinsilver}, but with many changes, not only to deal
with fixed jumps but also to clarify some obscure arguments in \cite{feinsilver}. To save the space,
we rely on some results proved in the first half of \cite{feinsilver} (which is relatively easier than the second half), and
have to omit some tedious computations after having stated the main technical points.

All processes are assumed to be defined on the infinite time interval $\bR_+=[0,\,\infty)$, but it is clear that the results also hold on a finite time interval.
For a manifold $X$, let $\myB(X)$ be the Borel $\sigma$-field on $X$ and let $\myB_+(X)$ be the space of nonnegative Borel functions on $X$.
Let $C(X)$, $C_b(X)$ and $C_c^\infty(X)$ be respectively the spaces of continuous functions, bounded continuous functions, and smooth functions
with compact supports on $X$.

\section{Inhomogeneous \levy processes in Lie groups} \label{seclp}

Let $x_t$ be an inhomogeneous \levy process in a Lie group $G$. By definition, $x_t$ is a rcll process in $G$ with independent
increments, that is, $x_s^{-1}x_t$ is independent of $\myF_s^x$ for $s<t$. It becomes a \levy process in $G$ if it also has stationary
increments, that is, if the distribution of $x_s^{-1}x_t$ depends only on $t-s$. Let $\mu_{s,t}$ be the distribution of $x_s^{-1}x_t$.
Then for $f\in\myB_+(G)$,
\[E[f(x_t)\mid\myF_s^x] = E[f(x_sx_s^{-1}x_t)\mid\myF_s^x] = \int_G f(x_sy)\mu_{s,t}(dy).\]
This shows that $x_t$ is an inhomogeneous Markov process with transition function $P_{s,t}$ given by $P_{s,t}f(x)=\int f(xy)\mu_{s,t}(dy)$.

Note that $x_t$ is left invariant in the sense that its transition function $P_{s,t}$ is left invariant,
that is, $P_{s,t}(f\circ l_g)=(P_{s,t}f)\circ l_g$ for $f\in\myB_+(G)$ and $g\in G$, where $l_g$ is the left translation $x\mapsto gx$
on $G$. Conversely, if $x_t$ is a left invariant inhomogeneous Markov process in $G$, then $E[f(x_s^{-1}x_t)\mid\myF_s^x]=P_{s,t}(f\circ l_{x_s^{-1}})
(x_s)=P_{s,t}f(e)$, where $e$ is the identity element of $G$. This implies that $x_t$ has independent increments.
We have proved the following result.

\begin{prop} \label{prinhlevyG}
A rcll process $x_t$ in $G$ is an inhomogeneous \levy process if and only
if it is a left invariant inhomogeneous Markov process.
\end{prop}

Note that the proof of Proposition~\ref{prinhlevyG} may be slightly modified to show that a rcll process in $G$ is a \levy process if
and only it is a Markov process with a left invariant transition function $P_t$ (see
also \cite[Proposition~1.2]{Liao2004}).

A measure function on $G$ is a family of $\sigma$-finite measures $\eta(t,\cdot)$
on $G$, $t\in\bR_+$, such that $\eta(s,\cdot)\leq\eta(t,\cdot)$ for $s<t$,
and $\eta(t,\cdot)\downarrow \eta(s,\cdot)$ as $t\downarrow s\geq 0$.
\index{measure function} Here the limit is set-wise, that is, $\eta(t,B)
\to\eta(s,B)$ for $B\in\myB(G)$.
The left limit $\eta(t-,\cdot)$ at $t>0$, defined as the
nondecreasing limit of measures $\eta(s,\cdot)$ as $s\uparrow t$,
exists and is $\leq\eta(t,\cdot)$.

A measure function $\eta(t,\cdot)$ may be regarded as a $\sigma$-finite measure on $\bR_+\times G$
and may be written as $\eta(dt,dx)$, given by $\eta((s,\,t]\times B)=\eta(t,B)- \eta(s,B)$
for $s<t$ and $B\in\myB(G)$. Conversely, any measure $\eta$ on $\bR_+\times G$ such
that $\eta(t,\cdot)=\eta([0,\,t]\times\cdot)$ is a $\sigma$-finite measure on $G$ for any $t>0$ may be identified
with the measure function $\eta(t,\cdot)$.

A measure function $\eta(t,\cdot)$ is called continuous at $t>0$ if $\eta(t,\cdot)=\eta(t-,\cdot)$,
and continuous if it is continuous at all $t>0$. In general, the set $J=\{t>0$: $\eta(t,G)>
\eta(t-,G)\}$, of discontinuity times, is at most countable, and
\begin{equation}
\eta(t,\cdot) = \eta^c(t,\cdot) + \sum_{s\leq t,\,s\in
J}\eta(\{s\}\times\cdot), \label{etaetaprimeeta}
\end{equation}
where $\eta^c(t,\cdot)=\int_{[0,\,t]\cap J^c}\eta(ds,\cdot)$ is a continuous measure function, called
the continuous part of $\eta(t,\cdot)$, and $\eta(\{s\}\times\cdot)=\eta(s,\cdot)-\eta(s-,\cdot)$

Recall $e$ is the identity element of $G$. The jump intensity measure
of an inhomogeneous \levy process $x_t$ is the measure function $\eta(t,\cdot)$
on $G$ defined by
\begin{equation}
\eta(t,B) = E\{\#\{s\in (0,\,t];\ \ x_{s-}^{-1}x_s\in B\ {\rm and}\
x_{s-}^{-1}x_s\neq e\}\},\ \ \ \ B\in\myB(G), \label{etatB}
\end{equation}
the expected number of jumps in $B$ by time $t$. The required $\sigma$-finiteness
of $\eta(t,\cdot)$ will be clear from Proposition~\ref{pretantoeta} later,
and then the required right continuity, $\eta(t,\cdot)\downarrow\eta(s,\cdot)$ as $t\downarrow s\geq 0$,
follows from (\ref{etatB}). It is clear that the process $x_t$ is continuous if and only if $\eta=0$.

Note that $\eta(0,\cdot)=0$ and $\eta(t,\{e\})=0$, and for $t>0$, $\eta(\{t\}\times G)\leq 1$ and
\begin{equation}
\nu_t = \eta(\{t\}\times\cdot) + [1-\eta(\{t\}\times G)]\delta_e\ \ \ \ \mbox{(where $\delta_e$
is the unit mass at $e$)} \label{nut}
\end{equation}
is the distribution of $x_{t-}^{-1}x_t$, so $\eta(t,\cdot)$ is continuous if
and only if $x_t$ is stochastically continuous.

Let $\{\xi_1,\ldots,\xi_d\}$ be a basis of the Lie algebra $\myg$ of $G$.
We will write $e^\xi$ for the exponential map $\exp(\xi)$ on $G$.
There are $\phi_1,\ldots,\phi_d\in C_c^\infty(G)$ such that $x=
e^{\sum_{i=1}^d\phi_i(x)\xi_i}$ for $x$ near $e$, called coordinate
functions associated to the basis $\{\xi_i\}$ of $\myg$. Note that $\xi_i\phi_j(e)=\delta_{ij}$.
The $\phi_\cdot$-truncated mean, or simply the mean, of a $G$-valued
random variable $x$ or its distribution $\mu$ is defined to be
\begin{equation}
b = e^{\sum_{j=1}^d\mu(\phi_j)\xi_j}\ \ \ \ \mbox{(where $\mu(\phi_j)=\int\phi_jd\mu$).} \label{meanb}
\end{equation}
The distribution $\mu$ is called small if its mean $b$ has coordinates $\mu(\phi_1),\ldots,
\mu(\phi_d)$, that is,
\begin{equation}
\phi_j(b) = \mu(\phi_j),\ \ \ \ 1\leq j\leq d. \label{phijb}
\end{equation}
This is the case when $\mu$ is sufficiently concentrated near $e$.

As defined in \cite{feinsilver}, a \levy measure function on $G$ is a continuous measure
function $\eta(t,\cdot)$ such that \vspace{2ex}

\noindent (a) \ $\eta(0,\cdot)=0$, $\eta(t,\{e\})=0$ and $\eta(t,U^c)<\infty$
for any $t\in\bR_+$ and neighborhood $U$ of $e$.

\noindent (b) \ $\eta(t,\Vert\phi_\cdot\Vert^2)<\infty$ for any $t\geq 0$,
where $\Vert\phi_\cdot\Vert^2=\sum_{i=1}^d\phi_i(x)^2$.
\vspace{2ex}

The notion of \levy measure functions is now extended. A measure
function $\eta(t,\cdot)$ on $G$ is called an extended \levy measure
function if (a) above and (b$'$) below hold. \vspace{2ex}

\noindent (b$'$) \ For $t\geq 0$, $\eta^c(t,\Vert\phi_\cdot\Vert^2)<\infty$, $\eta(\{t\}\times G)
\leq 1$ and with $\nu_t$ given by (\ref{nut}),
\begin{equation}
\sum_{s\leq t}\nu_s(\Vert\phi_\cdot-\phi_\cdot(h_s)\Vert^2) \ <\ \infty, \label{sumnu}
\end{equation}
where $\eta^c(t,\cdot)$ is the continuous part of $\eta(t,\cdot)$ and $h_t=\exp[\sum_j\nu_t(\phi_j)\xi_j]$ is the mean of $\nu_t$.
\vspace{2ex}

Note that $\nu_t=\delta_e$ and $h_t=e$ at a continuous point $t$ of $\eta(t,\cdot)$,
and hence the sum $\sum_{s\leq t}$ in (\ref{sumnu}) has at most countably many nonzero terms.
If $\eta(t,\cdot)$ is continuous, then (b$'$) becomes (b), and hence
a continuous extended \levy measure function is a \levy measure
function.

It can be shown directly that conditions (b) and (b$'$) are independent of the choice for coordinate functions $\phi_i$ and
basis $\{\xi_i\}$. This is also a consequence of Theorem~\ref{thmartrep} below.

A continuous path $b_t$ in $G$ with $b_0=e$ will be called a drift. A $d\times d$ symmetric matrix-valued function $A(t)=\{A_{jk}(t)\}$ will
be called a covariance matrix function if $A(0)=0$, $A(t)-A(s)\geq 0$ (nonnegative definite) for $s<t$, and $t\mapsto A_{jk}(t)$ is continuous.
A triple $(b,A,\eta)$ of a drift $b_t$, a covariance matrix function $A(t)$ and a \levy measure function $\eta(t,\cdot)$ will be called a
\levy triple on $G$. For $g\in G$, the adjoint map $\Ad(g)$: $\myg\to\myg$ is the differential of the conjugation $x\mapsto gxg^{-1}$ at $x=e$.
It is shown in \cite{feinsilver} that if $x_t$ is a stochastically continuous
inhomogeneous \levy process in $G$ with $x_0=e$, then
there is a unique \levy triple $(b,A,\eta)$ such that $x_t=z_tb_t$ and
\begin{eqnarray}
\forall f\in C_c^\infty(G), && f(z_t) -
\frac{1}{2}\int_0^t\sum_{j,k=1}^d[\Ad(b_s)\xi_j][\Ad(b_s)\xi_k]
f(z_s)\,dA_{jk}(s) \nonumber \\
&-& \int_0^t\int_G\{f(z_sb_s\tau b_s^{-1})-f(z_s)-\sum_{j=1}^d\phi_j(\tau)
[\Ad(b_s)\xi_j]f(z_s)\}\eta(ds,d\tau) \label{mart}
\end{eqnarray}
is a martingale under $\myF_t^x$. Conversely, given a \levy triple $(b,A,\eta)$,
there is a stochastically continuous inhomogeneous \levy
process $x_t$ with $x_0=e$ represented as above, unique in distribution.

Therefore, a stochastically continuous inhomogeneous \levy process in $G$ is represented
by a triple $(b,A,\eta)$ just like its counterpart in $\bR^d$. The complicated form of
the martingale in (\ref{mart}) with the presence of the drift $b_t$, as compared with
its counterpart (\ref{martRn}) on $\bR^d$, is caused by the
non-commutativity of $G$. This representation is extended to all
inhomogeneous \levy processes in $G$, not necessarily stochastically continuous,
in Theorem~\ref{thmartrep} below.

An extended drift on $G$ is a rcll path $b_t$ in $G$ with $b_0=e$. A
triple $(b,A,\eta)$ of an extended drift $b_t$, a covariance matrix
function $A(t)$ and an extended \levy measure function $\eta(t,\cdot)$ will be called
an extended \levy triple on $G$ if $b_{t-}^{-1}b_t=h_t$ for any $t>0$.

We note that for $G=\bR^d$, our definition of an extended \levy triple corresponds
to the assumptions in \cite[Theorem~II\,5.2]{js}. In particular, (a) and
(b$'$) corresponds to (i) - (iii), and $b_{t-}^{-1}b_t=h_t$ to (v),
but (iv) in \cite{js} is redundant as it is implied by the other conditions.

A rcll process $x_t$ in $G$ is
said to be represented by an extended \levy triple $(b,A,\eta)$ if
with $x_t=z_tb_t$,
\begin{eqnarray}
M_tf &=& f(z_t) - \frac{1}{2}
\int_0^t\sum_{j,k=1}^d[\Ad(b_s)\xi_j][\Ad(b_s)\xi_k]
f(z_s)dA_{jk}(s) \nonumber \\
&&\ \ -
\int_0^t\int_G\{f(z_sb_sxb_s^{-1})-f(z_s)-\sum_{j=1}^d\phi_j(x)
[\Ad(b_s)\xi_j]f(z_s)\}\eta^c(ds,dx) \nonumber \\
&&\ \ - \sum_{u\leq
t}\int_G[f(z_{u-}b_{u-}xh_u^{-1}b_{u-}^{-1})-f(z_{u-})] \nu_u(dx)
\label{Mtf}
\end{eqnarray}
is a martingale under the natural filtration $\myF_t^x$ of $x_t$ for
any $f\in C_c^\infty(G)$.

We note that in (\ref{Mtf}), the $\eta^c$-integral is absolutely integrable
and the sum $\sum_{u\leq t}$ converges absolutely a.s., and hence $M_tf$ is
a bounded random variable. This may
be verified by (b$'$) and Taylor's expansions of $f(z_sb_sxb_s^{-1})=
f(z_sb_se^{\sum_j\phi_j(x)\xi_j}b_s^{-1})$ at $x=e$
and $f(z_{u-}b_{u-}xh_u^{-1}b_{u-}^{-1})$ at $x=h_u$ (see the computation
in the proof of Lemma~\ref{leztmart}).

\begin{theo} \label{thmartrep}
Let $x_t$ be an inhomogeneous \levy process in $G$ with $x_0=e$.
Then there is a unique extended \levy triple $(b,A,\eta)$ on $G$ such
that $x_t$ is represented by $(b,A,\eta)$ as defined above.
Moreover, $\eta(t,\cdot)$ is the jump intensity measure of process
$x_t$ given by (\ref{etatB}). Consequently, $x_t$ is stochastically continuous if
and only if $(b,A,\eta)$ is a \levy triple.

  Conversely, given an extended \levy triple $(b,A,\eta)$ on $G$, there is
an inhomogeneous \levy process $x_t$ in $G$ with $x_0=e$, unique in distribution,
that is represented by $(b,A,\eta)$.
\end{theo}

\noindent {\bf Remark 1} \ As the jump intensity measure, $\eta(t,\cdot)$
is clearly independent of the choice for the basis $\{\xi_j\}$ of $\myg$
and coordinate functions $\phi_j$. By Lemma~\ref{leAtntoAt}, $A(t)$ is independent
of $\{\phi_j\}$ and the operator $\sum_{j,k=1}^dA_{jk}(t)\xi_j\xi_k$ is independent
of $\{\xi_j\}$.
\vspace{2ex}

In Theorem~\ref{thmartrep}, the representation of process $x_t$ is
given in the form of a martingale property of the shifted process $z_t=x_tb_t^{-1}$.
By Theorem~\ref{thmartfinvar} below, when the drift $b_t$ has a finite variation,
a form of martingale property holds directly for $x_t$.

A rcll path $b_t$ in a manifold $X$ is said to have a finite
variation if for any $f\in C_c^\infty(X)$, $f(b_t)$ has a finite
variation on any finite $t$-interval. Let $\xi_1,\ldots,\xi_d$ be a
family of smooth vector fields on $X$ such that $\xi_1(x),\ldots,\xi_d(x)$
form a basis of the tangent space $T_xX$ of $X$ at any $x$.
If $b_t$ is a continuous path in $X$ with a
finite variation, then there are uniquely defined real valued continuous
functions $b_j(t)$ of finite variation, $1\leq j\leq d$, with
$b_j(0)=0$, such that
\begin{equation}
\forall f\in C_c^\infty(X),\ \ \ \ f(b_t) = f(b_0) +
\int_0^t\sum_j\xi_jf(b_s)db_j(s). \label{xifdb}
\end{equation}
Indeed, $\xi_j(x)=\sum_k\alpha_{jk}(x) (\partial/\partial\psi_k)$
under local coordinates $\psi_1,\ldots,\psi_d$ on $G$,
and $df(b_t)=\sum_j(\partial/\partial \psi_j)f(b_t)d\psi_j(b_t)=\sum_{j,k}
\beta_{jk}(b_t)\xi_kf(b_t)d\psi_j(b_t)$,
where $\{\beta_{jk}(x)\}= \{\alpha_{jk}(x)\}^{-1}$, then $b_j(t)$ are determined
by $db_j(t)= \sum_k\beta_{kj}(b_t)d\psi_k(b_t)$. Conversely,
given $b_j(t)$ as above, a continuous path $b_t$ of finite variation satisfying (\ref{xifdb})
may be obtained by solving the integral equation
\[\psi_j(b_t) - \psi_j(b_{t_0}) = \sum_k\int_{t_0}^t \alpha_{kj}(b_s)db_k(s)\]
for $\psi_j(b_t)$ by the usual successive approximation method.

The functions $b_j(t)$ above will be called components of the path $b_t$ under the vector fields $\xi_1,\ldots,\xi_d$. When $X=G$,
these vector fields will be the basis of $\myg$ chosen before.

More generally, if $b_t$ is a rcll path in $X$ of finite variation,
then there is a unique continuous path $b_t^c$ in $X$ of finite
variation with $b_0^c=b_0$ such that, letting $b_j(t)$ be
 the components of $b_t^c$,
\begin{equation}
\forall f\in C_c^\infty(X),\ \ \ \ f(b_t) = f(b_0) +
\int_0^t\sum_j\xi_jf(b_s)db_j(s) + \sum_{s\leq t} [f(b_s) -
f(b_{s-})]. \label{xifdb2}
\end{equation}
To prove this, cover the path by finitely many coordinate
neighborhoods and then prove the claim on a Euclidean space. The
path $b_t^c$ will be called the continuous part of $b_t$.

\begin{theo} \label{thmartfinvar}
Let $x_t$ be an inhomogeneous \levy process in $G$ with $x_0=e$, represented by an extended \levy triple $(b,A,\eta)$. Assume $b_t$ is of finite variation. Then
\begin{eqnarray}
&& \hspace{-0.2in} f(x_t) - \int_0^t\sum_j \xi_jf(x_s)db_j(s) - \int_0^t\frac{1}{2}\sum_{j,k}\xi_j\xi_k f(x_s)\,dA_{jk}(s) \nonumber \\
&& - \int_0^t\int_G\{f(x_s\tau)-f(x_s)-\sum_j\phi_j(\tau)\xi_jf(x_s)\}\eta^c(ds,d\tau) \nonumber \\
&& - \sum_{u\leq t}\int_G[f(x_{u-}\tau)-
f(x_{u-})]\nu_u(d\tau) \label{martfinvarb}
\end{eqnarray}
is a martingale under $\myF_t^x$ for any $f\in C_c^\infty(G)$.

  Conversely, given an extended \levy triple $(b,A,\eta)$
with $b_t$ of finite variation, there is an inhomogeneous \levy process $x_t$
in $G$ with $x_0=e$, unique in distribution, such that (\ref{martfinvarb})
is a martingale under $\myF_t^x$ for $f\in C_c^\infty(G)$.
\end{theo}

The above theorem follows directly from Theorem~\ref{thmartrep} and the next lemma. Note that because $b_t$ has a finite variation and $\phi_i(b_u^{-1}b_u)=
\nu_u(\phi_i)$ for all but finitely many $u\leq t$, $\sum_{u\leq t}\vert\nu_u(\phi_i)\vert<\infty$ and hence $\sum_{u\leq t}\vert\int_G[f(x_{u-}\tau) -
f(x_{u-})]\nu_u(d\tau)\vert$ $<\infty$. The absolute integrability of the $\eta^c$-integral in (\ref{martfinvarb}) can be verified as in (\ref{Mtf}).

\begin{lemm} \label{lemartfinvar}
Let $(b,A,\eta)$ be an extended \levy triple on $G$ with $b_t$ of
finite variation. For an inhomogeneous \levy process $x_t=z_tb_t$
in $G$ with $x_0=e$, (\ref{Mtf}) being a martingale under $\myF_t^x$
for all $f\in C_c^\infty(G)$ is equivalent to (\ref{martfinvarb}) being a martingale
under $\myF_t^x$ for all $f\in C_c^\infty(G)$.
\end{lemm}

\noindent {\bf Proof} \ Let us assume (\ref{Mtf}) is a martingale.
Let $\Delta_n$: $0=t_{n0}<t_{n1}<\cdots<t_{ni}\uparrow\infty$ (as $i\uparrow\infty$) be a sequence
of partitions of $\bR_+$
with mesh $\Vert\Delta_n\Vert=\sup_{i\geq 1}(t_{ni}-t_{n\,i-1})\to 0$ as $n\to\infty$, and
let $f_{ni}(z)=f(zb_{t_{ni}})$. Then
\[f_{ni}(z_t) = f_{ni}(e) + M_t^{ni} + \int_{[0,\,t]} L(s,ds)f_{ni},\]
where $M_t^{ni}$ is a martingale with $M_0^{ni}=0$ and
\begin{eqnarray*}
&& L(t,dt)f = \frac{1}{2}\sum_{j,k}[\Ad(b_t)\xi_j][\Ad(b_t)\xi_k]
f(z_t)dA_{jk}(t) + \int_G\{f(z_tb_t\tau b_t^{-1})-f(z_t) \\
&& - \sum_j\phi_j(\tau) [\Ad(b_t)\xi_j]f(z_t)\}\eta^c(dt,d\tau)  +
\sum_{s\in dt} \int_G[f(z_{s-}b_{s-}\tau h_s^{-1}b_{s-}^{-1}) -
f(z_{s-})]\nu_s(d\tau).
\end{eqnarray*}
Let $b_t^n=b_{t_{ni}}$ for $t\in [t_{ni},\,t_{n\,i+1})$. Let $J$ be the set of fixed jump times of $x_t$.
We may assume $J\subset\Delta_n$ as $n\to\infty$ in the sense that $\forall u\in J$, $u\in\Delta_n$
for large $n$. Then $b_t^n$ is a step function and $b_t^n\to b_t$ as $n\to\infty$
uniformly for bounded $t$. It follows that for $t\in [t_{ni},\,t_{n\,i+1})$,
\begin{eqnarray*}
&& f(z_tb_t^n) = f(z_tb_{t_{ni}}) = f(e) +
\sum_{j=1}^i[f(z_{t_{nj}}b_{t_{nj}})-
f(z_{t_{n\,j-1}}b_{t_{n\,j-1}})] + [f(z_tb_{t_{ni}})-f(z_{t_{ni}}b_{t_{ni}})] \\
&=& f(e) + \sum_{j=1}^i[f(z_{t_{nj}}b_{t_{n\,j-1}})-
f(z_{t_{n\,j-1}}b_{t_{n\,j-1}})] + [f(z_tb_{t_{ni}}) - f(z_{t_{ni}}b_{t_{ni}})] \\
&&\ \ + \sum_{j=1}^i[f(z_{t_{nj}}b_{t_{nj}})-f(z_{t_{nj}}b_{t_{n\,j-1}})] \\
&=& f(e) + \sum_{j=1}^i[M_{t_{nj}}^{n\,j-1}-M_{t_{n\,j-1}}^{n\,j-1}] +
  \sum_{j=1}^i\int_{(t_{n\,j-1},\,t_{nj}]}L(s,ds)f_{n\,j-1} + [M_t^{ni}-M_{t_{ni}}^{ni}] \\
&& \ \ + \int_{(t_{ni},\,t]}L(s,ds)f_{ni} + \sum_{j=1}^i
[f(z_{t_{nj}}b_{t_{nj}})-f(z_{t_{nj}}b_{t_{n\,j-1}})] \\
&=& f(e) + M_t^{(n)} + \int_{[0,\,t]}L(s,ds)(f\circ r_{b_{s-}^n}) +
\sum_{j=1}^i [f(z_{t_{nj}}b_{t_{nj}})-f(z_{t_{nj}}b_{t_{n\,j-1}})],
\end{eqnarray*}
where $M_t^{(n)}=\sum_{j=1}^i[M_{t_{nj}}^{n\,j-1}-M_{t_{n\,j-1}}^{n\,j-1}] +
[M_t^{ni}-M_{t_{ni}}^{ni}]$, $t\in[t_{ni},\,t_{n\,i+1})$, is a martingale, and $r_b$ is
the right translation $x\mapsto xb$ on $G$. As $n\to\infty$,
by the uniform convergence $b_t^n\to b_t$, $f(z_tb_t^n)\to
f(z_tb_t)=f(x_t)$, and by the continuity of $A(t)$ and $\eta^c(t,\cdot)$,
\begin{eqnarray*}
&& \int_{[0,\,t]}L(s,ds)(f\circ r_{b_{s-}^n}) \ \to\
\int_{[0,\,t]}L(s,ds) (f\circ r_{b_{s-}}) =
\frac{1}{2}\int_0^t\sum_{p,q}\xi_p\xi_qf(x_s)dA_{pq}(s) + \\
&& \int_0^t\int_G[f(x_s\tau) - f(x_s) - \sum_p\phi_p(\tau)
\xi_pf(x_s)] \eta^c(ds,d\tau) + \sum_{s\leq t}\int_G[f(x_{s-}\tau
h_s^{-1}) - f(x_{s-})]\nu_s(d\tau)
\end{eqnarray*}
and, by (\ref{xifdb2}) and noting that $b_p(t)$ is continuous in $t$ and $J\subset\Delta_n$ as $n\to\infty$,
\begin{eqnarray*}
&& \hspace{-0.5in} \sum_{j=1}^i[f(z_{t_{nj}}b_{t_{nj}})-f(z_{t_{nj}}
b_{t_{n\,j-1}})] = \sum_{j=1}^i\{\int_{t_{n\,j-1}}^{t_{nj}}\sum_p
\xi_pf(z_{t_{nj}}b_s)db_p(s) + \sum_{t_{n\,j-1}<u\leq t_{nj}}[f(z_{t_{nj}}b_u)-f(z_{t_{nj}}b_{u-})]\} \\
&\to& \int_0^t\sum_p\xi_pf(z_sb_s)db_p(s) + \sum_{u\leq t}[f(z_ub_u)-f(z_ub_{u-})] \\
&=& \int_0^t\sum_p\xi_pf(z_sb_s)db_p(s) + \sum_{u\leq t}[f(x_u) - f(x_uh_u^{-1})] \\
&=& \int_0^t\sum_p\xi_pf(x_s)db_p(s) + \sum_{u\leq t}\int_G[f(x_{u-}\tau)-f(x_{u-}\tau h_u^{-1})]\nu_u(d\tau) + M_t',
\end{eqnarray*}
where $M_t'=\sum_{u\leq t}\{f(x_u)-f(x_uh_u^{-1})-E[f(x_u)-f(x_uh_u^{-1})\mid\myF_{u-}^x]\}$ is a bounded martingale.
% (the boundedness follows from the finite variation of $b_t$ noting $h_t=b_{t-}^{-1}b_t$)
It follows that $M_t^{(n)}$ converges boundedly to a martingale $M_t$ given by (\ref{martfinvarb}).

Now we assume (\ref{martfinvarb}) is a martingale. Because $z_t=x_tb_t^{-1}$,
the above computation can be repeated with $x_t$
and $z_t$ switching roles, and $b_t$ replaced by $b_t^{-1}$, to show
that (\ref{Mtf}) is a martingale. Note that the components $\beta_j(t)$ of the continuous part
of $b_t^{-1}$ under the basis $\{\xi_j\}$ are given by
\begin{equation}
  d\beta_j(t) = -\sum_k[\Ad(b_t)]_{jk}db_k(t), \label{dbetai}
\end{equation}
where $[\Ad(b)]_{jk}$ is the matrix representing $\Ad(b)$, that is, $\Ad(b)\xi_k=\sum_j[\Ad(b)]_{jk}\xi_j$.
% , because,
% denoting $J$: $G\to G$ as the inverse map,
% \begin{eqnarray*}
% && f(b_t^{-1}) = (f\circ J)(b_t) = f(b_0^{-1}) + \int_0^t\sum_i
% \xi_i(f\circ J)(b_s)db_i(s) + \sum_{u\leq t}[f(b_u^{-1})-
% f(b_{u-}^{-1})] \\
% &=& f(b_0^{-1}) + \int_0^t\sum_i[-\Ad(b_t)\xi_i]f(b_s^{-1})db_i(s) +
% \sum_{u\leq t}[f(b_u^{-1})-f(b_{u-}^{-1})]. \ \ \Box
% \end{eqnarray*}
\ $\Box$

\section{Measure functions} \label{secmeasfun}

We now let $J=\{u_1,u_2,u_3,\ldots\}$ be the set of fixed jump times
of an inhomogeneous \levy process $x_t$ (the list may be finite or
empty, and may not be ordered in magnitude). As in the proof of Lemma~\ref{lemartfinvar},
let $\Delta_n$: $0=t_{n0}<t_{n1}<t_{n2}<\cdots< t_{nk}\uparrow\infty$ (as $k\uparrow\infty$)
be a sequence of partitions of $\bR_+$ with mesh $\Vert\Delta_n\Vert\to 0$ as $n\to\infty$, such
that $\Delta_n$ contains $J$ as $n\to\infty$ in the sense that any $u\in J$ is in $\Delta_n$
when $n$ is sufficiently large. For $m>0$, let $J_m=\{u_1,u_2,\ldots,u_m\}$.

For $i=1,2,3,\ldots$, let $x_{ni}=x_{t_{n\,i-1}}^{-1}x_{t_{ni}}$.
Then for each fixed $n$, $x_{ni}$ are independent random variables
in $G$. Let $\mu_{ni}$ be their distributions, and define the measure function
\begin{equation}
\eta_n(t,\cdot) = \sum_{t_{ni}\leq t}\mu_{ni}\ \ \ \ \mbox{(setting $\eta_n(0,\cdot)=0$)}. \label{etant}
\end{equation}

\begin{prop} \label{pretanbdd}
For any $T>0$ and any neighborhood $U$ of $e$, $\eta_n(T,U^c)$ is bounded in $n$ and $\eta_n(T,U^c)\downarrow 0$ as $U\uparrow G$ uniformly in $n$.
\end{prop}

\noindent {\bf Proof} \ We first establish an equi-continuity type property for $\eta_n(t,\cdot)$.

\begin{lemm} \label{leetanequicont}
For any $T>0$, neighborhood $U$ of $e$ and $\varepsilon>0$, there are integers $n_0$ and $m>0$, and $\delta>0$, such that if $n\geq n_0$
and $s,t\in\Delta_n\cap[0,\,T]$ with $0<t-s<\delta$ and $(s,\,t]\cap J_m=\emptyset$ (the empty set), then $\eta_n(t,U^c)-\eta_n(s,U^c)<\varepsilon$.
\end{lemm}

By Borel-Cantelli Lemma, the independent increments and rcll paths
imply that for any neighborhood $U$ of $e$,
\begin{equation}
\sum_{u\leq T,\,u\in J}P(x_{u-}^{-1}x_u\in U^c) < \infty.
\label{sumPxuj}
\end{equation}

  Suppose the claim of the lemma is not true. Then for some $\varepsilon>0$,
and for any $n_0,m$ and $\delta>0$, there are $s_n,t_n
\in\Delta_n\cap[0,\,T]$ with $n\geq n_0$, $t_n-s_n<\delta$ and
$(s_n,\,t_n]\cap J_m=\emptyset$ such that
$\eta_n(t_n,U^c)-\eta_n(s_n,U^c)\geq\varepsilon$.
Letting $n_0\to\infty$ and $\delta\to 0$ yields a subsequence of
$n\to\infty$ such that $s_n$ and $t_n$ converge to a common limit
$t\leq T$ as $n\to\infty$.

Let $s_n=t_{ni}$ and $t_n=t_{nj}$. In the following computation,
we will write $A_n\approx B_n$ if there is
a constant $c>0$ such that $(1/c)A_n\leq B_n\leq cA_n$. Then
\begin{eqnarray*}
&& \eta_n(t_n,U^c) - \eta_n(s_n,U^c) = \sum_{p=i+1}^jP(x_{t_{n\,p-1}}^{-1}x_{t_{np}}\in U^c)
\approx - \sum_{p=i+1}^j\log[1 - P(x_{t_{n\,p-1}}^{-1}x_{t_{np}}\in U^c)] \\
&=& -\log P[\cap_{p=i+1}^j(x_{t_{n\,p-1}}^{-1}x_{t_{np}}\in U)]
\ =\ -\log P(A_n),
\end{eqnarray*}
where $A_n=\cap_{p=i+1}^j(x_{t_{n\,p-1}}^{-1}x_{t_{np}}\in U)$. Because $\eta_n(t_n,U^c)-\eta_n(s_n,U^c)
\geq\varepsilon$, $P(A_n)\leq 1-\varepsilon_1$ for some constant $\varepsilon_1>0$.
Then $P(A_n^c)>\varepsilon_1$. Note that on $A_n^c$,
the process $x_t$ makes a $U^c$-oscillation during the time interval $[s_n,\,t_n]$.

There are three possible cases and we will reach a contradiction in
all these cases. Case one, there are infinitely many $s_n\downarrow
t$. This is impossible by the right continuity of paths at time $t$. Case two,
there are infinitely many $t_n\uparrow t$. This is impossible by
the existence of path left limit at time $t$. Case three, there are
infinitely many $s_n<t\leq t_n$. This implies $P(x_{t-}^{-1}x_t\in
U^c)\geq\varepsilon_1$. Then $t\in J$. Because $(s_n,\,t_n]\cap
J_m=\emptyset$, $t\not\in J_m$. By (\ref{sumPxuj}), $m$ may be chosen so that $\sum_{u\leq T,\,u\in
J-J_m}P(x_{u-}^{-1}x_u\in U^c)<\varepsilon_1$, which is impossible.
Lemma~\ref{leetanequicont} is now proved.

To prove Proposition~\ref{pretanbdd}, fix $\varepsilon>0$ and let $n_0,m,\delta$
be as in Lemma~\ref{leetanequicont}. It suffices to prove for
$n\geq n_0$. For each $n\geq n_0$, $J_m$ may be covered by no more
than $m$ sub-intervals of the form $[t_{n\,i-1},\,t_{ni}]$, and the
rest of the interval $[0,\,T]$ may be covered by finitely many, say
$p$, sub-intervals $[s,\,t]$ as in Lemma~\ref{leetanequicont}. Then
$\eta_n(T,U^c)\leq p\varepsilon+m$, and hence $\eta_n(t,U^c)$
is bounded.

Now let $V$ be a neighborhood of $e$ such that $V^{-1}V=\{x^{-1}y$;
$x,y\in V\}\subset U$. Then
\begin{eqnarray*}
&& \eta_n(T,U^c) = \sum_{t_{ni}\leq T} P(x_{t_{n\,i-1}}^{-1}x_{t_{ni}}\in U^c) \approx
-\sum_{t_{ni}\leq T}\log[1 - P(x_{t_{n\,i-1}}^{-1}x_{t_{ni}}\in U^c)] \\
&=& -\log P\{\cap_{t_{ni}\leq T}[x_{t_{n\,i-1}}^{-1}x_{t_{ni}}\in U]\}\ \leq\ -\log
P\{\cap_{t_{ni}\leq T}[x_{t_{n\,i-1}}\in V\ {\rm and}\ x_{t_{ni}}\in V]\} \\
&\leq& - \log P\{\cap_{t\leq T}[x_t\in V]\} \ \downarrow\ 0
\end{eqnarray*}
as $V\uparrow G$ uniformly in $n$ because $P\{\cap_{t\leq T}[x_t\in
V]\} \uparrow 1$. \ $\Box$ \vspace{2ex}

Let $\eta(t,\cdot)$ be the jump intensity measure of process $x_t$ defined by (\ref{etatB}).

\begin{prop} \label{pretantoeta}
For any $t>0$ and $f\in C_b(G)$ vanishing in a neighborhood of $e$,
\[\eta_n(t,f)\ \to\ \eta(t,f).\]
Moreover, $\eta(t,U^c)<\infty$ for any neighborhood $U$ of $e$.
\end{prop}

\noindent {\bf Proof} \ We may assume $f\geq 0$. Let $F=\sum_{u\leq t}f(x_{u-}^{-1}x_u)$ and $F_n=\sum_{i=1}^p f(x_{ni})$, where $p$ is
the largest index $i$ such that $t_{ni}\leq t$. Then $\eta(t,f)=E(F)$ and $\eta_n(t,f)=E(F_n)$.
Because $\Delta_n$ contains $J$ as $n\to\infty$, $F_n\to F$ a.s. as $n\to\infty$. By the independence of $x_{n1},x_{n2},x_{n3},\ldots$,
\begin{eqnarray*}
&& E(F_n^2) = E\{[\sum_{i=1}^p f(x_{ni})]^2\} = \sum_{i=1}^p
\mu_{ni}(f^2) + \sum_{i\neq j}\mu_{ni}(f)\mu_{nj}(f) \\
&=& \sum_{i=1}^p\mu_{ni}(f^2) + [\sum_{i=1}^p\mu_{ni}(f)]^2 -
\sum_{i=1}^p\mu_{ni}(f)^2 \leq \eta_n(t,f^2) + [\eta_n(t,f)]^2,
\end{eqnarray*}
which is bounded by Proposition~\ref{pretanbdd}. Therefore, $E(F_n^2)$
are uniformly bounded in $n$, and hence $F_n$ are
uniformly integrable. It follows that $\eta_n(t,f)=E(F_n)\to
E(F)=\eta(t,f)$. The finiteness of $\eta(t,U^c)$ now follows from
Proposition~\ref{pretanbdd}. \ $\Box$ \vspace{2ex}

  Define
\begin{equation}
\eta^m(t,B) = E[\#\{s\leq t;\ \ x_{s-}^{-1}x_s\in B,\
x_{s-}^{-1}x_s\neq e \ {\rm and}\ s\not\in J_m\}]. \label{etam}
\end{equation}
Then $\eta(t,\cdot)\geq\eta^m(t,\cdot)\downarrow\eta^c(t,\cdot)$ as $m\uparrow\infty$.
Let $\Delta_n^m$ be the subset of $\Delta_n$ consisting
of $t_{ni}$ such that $u\in (t_{n\,i-1},\,t_{ni}]$ for some $u\in J_m$, and
let $\eta_n^m(t,\cdot) = \sum_{t_{ni}\leq t,\,t_{ni}\not\in\Delta_n^m}\mu_{ni}$.

\begin{lemm} \label{leetanmto}
Fix $m>0$. Then as $n\to\infty$, $\eta_n^m(t,f)\to\eta^m(t,f)$ for
any $t>0$ and $f\in C_b(G)$ vanishing in a neighborhood of $e$.
\end{lemm}

\noindent {\bf Proof} \ Let $F'=\sum_{s\leq t,\,s\not\in J_m}f(x_{s-}^{-1} x_s)$ and $F_n'=\sum_{t_{ni}\leq t,\,t_{ni}\not\in\Delta_n^m}f(x_{ni})$.
Then $\eta^m(t,f)=E(F')$ and $\eta_n^m(t,f)=E(F_n')$. Because $\Delta_n$ contains $J$ as $n\to\infty$, $F_n'\to F'$ a.s. as $n\to\infty$.
One may assume $f\geq 0$. Because $F_n'\leq F_n$ in the notation of the proof of Proposition~\ref{pretantoeta}, so $F_n'$ are uniformly integrable
and hence $\eta_n^m(t,f)=E(F_n')\to E(F')=\eta^m(t,F)$. \ $\Box$
\vspace{2ex}

The notion of measure functions is now extended to matrix-valued measures. Thus a family of $d\times d$ symmetric matrix valued functions $A(t,B)=
\{A_{jk}(t,B)\}_{j,k=1,2,\ldots,d}$, for $t\in\bR_+$ and $B\in\myB(G)$, is called a matrix-valued measure function on $G$ if $A_{jk}(t,\cdot)$ is
a finite signed measure on $G$, $A(0,B)$ and $A(t,B)-A(s,B)$ are nonnegative definte for $s<t$, and $A_{jk}(t,B)\to A_{jk}(s,B)$ as $t\downarrow s$.

The trace of a matrix-valued measure function $A(t,\cdot)$, $q(t,\cdot)=\trace[A(t,\cdot)]$, is a finite measure function such
that for $s<t$ and $j,k=1,2,\ldots,d$,
\begin{equation}
\vert A_{jk}(t,\cdot)-A_{jk}(s,\cdot)\vert\ \leq\ q(t,\cdot)-
q(s,\cdot)\ =\ q((s,\,t]\times\cdot). \label{Aleqq}
\end{equation}

Let $A^n(t,\cdot)$ be the matrix-valued measure functions on $G$ defined by
\begin{equation}
A_{jk}^n(t,B) = \sum_{t_{ni}\leq t}\int_B[\phi_j(x)-\phi_j(b_{ni})]
[\phi_k(x)-\phi_k(b_{ni})]\mu_{ni}(dx), \label{AjkntB}
\end{equation}
where $b_{ni}=e^{\sum_j\mu_{ni}(\phi_j) \xi_j}$ is the mean of
$x_{ni}$. Its trace is
\begin{equation}
q^n(t,B) = \sum_{t_{ni}\leq t}\int_B\Vert \phi_\cdot(x)-
\phi_\cdot(b_{ni})\Vert^2\mu_{ni}(dx), \label{qntB}
\end{equation}
where $\phi_\cdot(x)=(\phi_1(x),\ldots,\phi_d(x))$ and $\Vert\cdot\Vert$
is the Euclidean norm on $\bR^d$.

In the rest of this section, let $\eta(t,\cdot)$ be the jump
intensity measure of an inhomogeneous \levy process $x_t$ in $G$, defined by (\ref{etatB}),
with $J$ being the set of its discontinuity
times. Let $\nu_t$ be defined by (\ref{nut}) with mean $h_t$. Then $\nu_t$ and $h_t$ are
respectively the distribution and the mean of $x_{t-}^{-1}x_t$,
which are nontrivial only when $t\in J$.
We will see that $\eta(t,\cdot)$ is an extended \levy measure function. This
fact is now proved below under an extra assumption.

\begin{prop} \label{pretaextlevy}
Assume $q^n(t,G)$ is bounded in $n$ for each $t>0$. Then the jump
intensity measure $\eta(t,\cdot)$ is an extended \levy measure
function, and hence its continuous part $\eta^c(t,\cdot)$ is a \levy
measure function.
\end{prop}

\noindent {\bf Proof} \ Let $c_q$ be a constant such that $q^n(t,G)\leq c_q$
for all $n$. Any $u\in J$ is contained in $(t_{n\,i-1},\,t_{ni}]$ for some $i=i_n$
and $\nu_u(\Vert\phi_\cdot-\phi_\cdot(h_u)\Vert^2)=\lim_{n\to\infty}
\mu_{ni}(\Vert\phi_\cdot-\phi_\cdot(b_{ni})\Vert^2)$. It follows
that $\sum_{u\leq t,u\in
J}\nu_u(\Vert\phi_\cdot-\phi_\cdot(h_u)\Vert^2)\leq c_q$. It remains
to prove the finiteness of $\eta^c(t,\Vert\phi_\cdot\Vert^2)$. Let
\begin{equation}
r_n^m = \max_{t_{ni}\leq t,\,t_{ni}\not\in
\Delta_n^m}\Vert\phi_\cdot(b_{ni})\Vert, \label{rmn}
\end{equation}
where $\Delta_n^m$ is defined after (\ref{etam}). Then $r_n^m=
\max_{t_{ni}\leq t,\,t_{ni}\not\in\Delta_n^m}\Vert\mu_{ni}(\phi_\cdot)\Vert$
for large $m$, and $\lims_{n\to\infty}r_n^m=\sup_{u\leq t,\,u\not\in J_m}
\Vert\nu_u(\phi_\cdot)\Vert\to 0$ as $m\to\infty$.
For two neighborhoods $U\subset V$ of $e$,
\begin{eqnarray*}
&& \eta^c(t,\Vert\phi_\cdot\Vert^21_{V^c}) \leq
\eta^m(t,\Vert\phi_\cdot\Vert^21_{V^c}) \ \leq\
\lims_{n\to\infty}\int_{U^c}\Vert\phi_\cdot(x)\Vert^2\eta_n^m(t,dx)\
\ \ \ \mbox{(by
Lemma~\ref{leetanmto})} \\
&\leq& 2\lims_{n\to\infty}\sum_{t_{ni}\leq t,\,t_{ni}\not\in
J_m}\int_{U^c}
\Vert\phi_\cdot(x)-\phi_\cdot(b_{ni})\Vert^2\mu_{ni}(dx) +
2\lims_{n\to\infty}(r_n^m)^2\eta_n(t,U^c) \\
&\leq& 2c_q +
2(\lims_{n\to\infty}(r_n^m)^2)(\lims_{n\to\infty}\eta_n(t,U^c)).
\end{eqnarray*}
Now letting $m\to\infty$ and then $V\downarrow\{e\}$ shows $\eta^c(t,\Vert\phi_\cdot\Vert^2)
\leq 2c_q$. \ $\Box$
\vspace{2ex}

Now consider the following equi-continuity type condition on $q^n(t,\cdot)$.
\vspace{2ex}

\noindent (A) \ For any $T>0$ and $\varepsilon>0$, there are integers $n_0,m>0$
and $\delta>0$ such that if $n\geq n_0$ and $s,t\in [0,\,T]$ with $0<t-s<\delta$
and $(s,\,t]\cap J_m=\emptyset$, then $q^n(t,G)- q^n(s,G)<\varepsilon$.

\begin{lemm} \label{leqnTGbdd}
Under (A), $q^n(T,G)$ is bounded in $n$, and $q^n(T,U^c)\downarrow 0$ uniformly in $n$ as $U\uparrow G$,
where $U$ is a neighborhood of $e$.
\end{lemm}

\noindent {\bf Proof} \ The boundedness of $q^n(T,G)$ is derived
from (A) in the same way as the boundedness of $\eta_n(T,U^c)$ in
Proposition~\ref{pretanbdd} is derived from Lemma~\ref{leetanequicont}.
Because the convergence $\eta_n(T,U^c)\downarrow 0$, $U\uparrow G$, is uniform
in $n$, so is $q^n(T,U^c)\downarrow 0$. \ $\Box$ \vspace{2ex}

\begin{lemm} \label{leAtntoAt}
Assume (A). Then there is a matrix valued measure function $A(t,\cdot)$
on $G$ such that along a subsequence of $n\to\infty$, $A^n(t,f)\to A(t,f)$
for all $t>0$ and $f\in C_b(G)$. Moreover, there is a covariance matrix function $A(t)$
on $G$ such that
\begin{eqnarray}
A_{jk}(t,f) &=& f(e)A_{jk}(t) + \int_Gf(x)\phi_j(x)\phi_k(x)
\eta^c(t,dx) \nonumber \\
&&\ \ + \sum_{u\leq t,\,u\in
J}\int_Gf(x)[\phi_j(x)-\phi_j(h_u)][\phi_k(x)-\phi_k(h_u)]\nu_u(dx).
\label{AtUAtintU}
\end{eqnarray}
Furthermore, $A(t)$ is independent of the choice for coordinate functions $\phi_j$
and the operator $\sum_{j,k=1}^d A_{jk}(t)\xi_j\xi_k$ is independent of the choice
for the basis $\{\xi_j\}$ of $\myg$.
\end{lemm}

\noindent {\bf Proof} \ Let $\Lambda$ be a countable dense subset of $[0,\,T]$ containing $J\cap [0,\,T]$
and let $H$ be a countable subset of $C_b(G)$. By the boundedness of $q^n(T,G)$ and (\ref{Aleqq}), along a subsequence of $n\to\infty$, $A^n(t,f)$ converges
for all $t\in \Lambda$ and $f\in H$. Let $K_n$ be an increasing sequence of compact subsets of $G$ such that $K_n\uparrow G$.
For each $K_n$, there is a countable subset $H_n$ of $C_b(K_n)$ that is dense in $C(K_n)$ under the supremum norm.
We will extend functions in $H_n$ to be functions on $G$ without increasing their suprenorms and let $H=\cup_nH_n$.

Because for any relatively compact neighborhood $U$ of $e$, $q^n(T,U^c)\downarrow 0$ as $U\uparrow G$ uniformly in $n$, it follows that $A^n(t,f)$ converges for
any $f\in C_b(G)$. By (A), $A^n(t,f)$ converges to some $A(t,f)$ for any $t\in [0,\,T]$ and $f\in C_b(G)$, and the convergence is uniform in $t$ and $f$ bounded
by a fixed constant. Because $A^n(t,f)$ is right continuous with left limits in $t$, so is $A(t,f)$
and hence $A(t,\cdot)$ is a matrix-valued measure function. Moreover, the jumps of $A(t,f)$
are $\int_Gf(x)[\phi_j(x)-\phi_j(h_t)][\phi_k(x)-\phi_k(h_t)]\nu_t(dx)$ at $t\in J$.

Let the sum $\sum_{t_{ni}\leq t}$ in (\ref{AjkntB}), which defines $A^n(t,\cdot)$, be broken into
two partial sums: $\sum_{t_{ni}\leq t,\,t_{ni}\not\in\Delta_n^m}$ and $\sum_{t_{ni}\leq t,
\,t_{ni}\in\Delta_n^m}$, and write $A^n(t,\cdot)=A^{n,m}(t,\cdot)+B^{n,m}(t,\cdot)$,
where $A^{n,m}(t,\cdot)=\sum_{t_{ni}\leq t,\,t_{ni}\not\in\Delta_n^m}$ and $B^{n,m}(t,\cdot)=\sum_{t_{ni}\leq t,
\,t_{ni}\in\Delta_n^m}$. Then for $f\in C_b(G)$, as $n\to\infty$,
\begin{equation}
B_{jk}^{n,m}(t,f)\ \to\ \sum_{s\leq t,\,s\in J_m}\int_Gf(x)[\phi_j(x)-\phi_j(h_t)][\phi_k(x)-\phi_k(h_t)]\nu_s(dx), \label{Bnmtfto}
\end{equation}
and hence $\lim_{n\to\infty}A^{n,m}(t,f)$ also exists. Note that the sum $\sum_{s\leq t,\,s\in J_m}$ in (\ref{Bnmtfto}) contains
the jumps of $A(s,f)$ at $s\in [0,\,t]\cap J_m$ and it converges to
\[B_{jk}'(t,f) = \sum_{s\leq t,\,s\in J}\int_Gf(x)[\phi_j(x)-\phi_j(h_t)][\phi_k(x)-\phi_k(h_t)]\nu_s(dx)\]
  as $m\to\infty$, and $B'(t,\cdot)$ is a matrix valued measure
function. It follows that as $m\to\infty$, $\lim_{n\to\infty}A^{n,m}(t,f)$ converges to a matrix valued measure
function $A'(t,f)=A(t,f)-B'(t,f)$.

Let $\psi\in C_c^\infty(G)$ with $\psi=1$ near $e$ and $0\leq\psi\leq 1$ on $G$.
Because $\eta_n(t,1-\psi)$ is bounded in $n$, and $\lims_{n\to\infty}r_n^m\to 0$
as $m\to\infty$ for $r_n^m$ defined by (\ref{rmn}), by Lemma~\ref{leetanmto},
\begin{equation}
A_{jk}^{n,m}(t,(1-\psi)f) = \int[1-\psi(x)]
f(x)\phi_j(x)\phi_k(x)\eta^c(t,dx) + r_{nm}' \label{Ajknmt}
\end{equation}
with $\lims_{n\to\infty}r_{nm}'\to 0$ as $m\to\infty$. Then
\[A_{jk}'(t,(1-\psi)f) = \int[1-\psi(x)]f(x)\phi_j(x)\phi_k(x)\eta^c(t,dx).\]
Let $\psi=\psi_p\downarrow 1_{\{e\}}$ with supp$(\psi_p)\downarrow\{e\}$ as $p\uparrow\infty$, and define $A(t)=
\lim_{p\to\infty}A'(t,\psi_p)$. Then $\lim_{p\to\infty}A'(t,f\psi_p)=f(e)A(t)$.
Because $A(t,f)=A'(t,f)+B'(t,f)$ and $A'(t,f)=A'(t,\psi_pf)+A'(t,(1-\psi_p)f))$,
letting $p\to\infty$ yields (\ref{AtUAtintU}). Because the jumps of $A(t,f)$
are accounted for by the sum $\sum_{u\leq t,\,u\in J}$ in (\ref{AtUAtintU}), it follows
that $A(t)$ is continuous and hence is a covariance matrix function.

Let $\{\tilde{\phi}_j\}$ be another set of coordinate functions associated to
the same basis $\{\xi_j\}$ of $\myg$. Then $\tilde{\phi}_j=\phi_j$ in a neighborhood $V$ of $e$.
Let $\tilde{b}_{ni}$, $\tilde{A}^{n,m}(t,\cdot)$ and $\tilde{A}(t)$ be $b_{ni}$, $A^{n,m}(t,\cdot)$
and $A(t)$ for $\tilde{\phi}_i$. Because for $t_{ni}\leq t$
with $t_{ni}\not\in\Delta_n^m$ and a large $m$, $\mu_{ni}$ is small in the sense defined
before (\ref{phijb}), $\phi_j(b_{ni})=\mu_{ni}(\phi_j)$
and $\tilde{\phi}_j(\tilde{b}_{ni})=\mu_{ni}(\tilde{\phi}_j)$. Then for $f\in C_b(G)$ vanishing on $V^c$,
\begin{eqnarray*}
\tilde{A}_{j,k}^{n,m}(t,f) &=& \sum_{t_{ni}\leq t,\,t_{ni}\not\in\Delta_n^m}\int_G\int_G\int_G f(x)
[\tilde{\phi}_j(x)-\tilde{\phi}_j(y)][\tilde{\phi}_k(x)-\tilde{\phi}_k(z)]\mu_{ni}(dx)\mu_{ni}(dy)\mu_{ni}(dz) \\
&=& \sum \int_V\int_V\int_Vf(x)[\phi_j(x)-\phi_j(y)][\phi_k(x)-\phi_k(z)]\mu_{ni}(dx)\mu_{ni}(dy)\mu_{ni}(dz) + R,
\end{eqnarray*}
where $R=\sum[\int_V\int_{V^c}\int_V+\int_{V^c}\int_V\int_V+\int_{V^c}\int_{V^c}\int_V]$, and a similar expression
holds for $A_{jk}^{n,m}(t,f)$. Subtract the two expressions, it can be shown that $\vert\tilde{A}_{jk}^{n,m}(t,f)-
A_{jk}^{n,m}(t,f)\vert$ is controlled by $\max_{t_{ni}\leq t,\,t_{ni}\not\in\Delta_n^m}[\mu_{ni}(V^c)+
\int\int\mu_{ni}(dx)\mu_{ni}(dy)\Vert\phi_\cdot(x)-\phi_\cdot(y)\Vert]\eta_n(t,V^c)$, and it follows
that $\lim_{m\to\infty}\lim_{n\to\infty}\vert\tilde{A}_{jk}^{n,m}(t,f)-A_{jk}^{n,m}(t,f)
\vert=0$. This implies that $\tilde{A}(t)=A(t)$.

Now let $\{\tilde{\xi}_j\}$ be another basis of $\myg$ such that $\xi_j=\sum_{k=1}^da_{jk}\tilde{\xi}_k$.
Then $\tilde{\phi}_j=\sum_{k=1}^da_{kj}\phi_k$ are the coordinate
functions associated to $\{\tilde{\xi}_j\}$, and from the above displayed expression for $\tilde{A}^{n,m}(t,f)$
in terms of $\tilde{\phi}_j$, $\tilde{A}_{jk}^{n,m}(t,f)=\sum_{p,q}a_{pj}a_{qk}A_{pq}^{n,m}(t,f)$.
This implies that $\tilde{A}_{jk}(t)=\sum_{p,q}a_{pj}a_{qk}A_{pq}(t)$,
and hence $\sum_{j,k}\tilde{A}_{jk}(t)\tilde{\xi}_j\tilde{\xi}_k=\sum_{j,k}A_{jk}(t)\xi_j\xi_k$. \ $\Box$
\vspace{2ex}

Let $Y$ be a smooth manifold equipped with a compatible metric $r$, which will be taken to be $G\times G$ in the proof of Lemma~\ref{leztmart}.
Let $y^n$ and $y$ be rcll functions: $\bR_+\to Y$. Assume for any $t>0$, $r(y^n(t_{ni}),y(t_{ni}))\to 0$ as $n\to\infty$ uniformly for $t_{ni}\leq t$.
Let $F(y,b,x)$ and $F_{jk}(y,b,x)$ be bounded continuous functions on $Y\times G\times G$.

\begin{lemm} \label{lefjkntofjk}
Assume the above and (A), and let $A(t)$ be the covariance matrix function in Lemma~\ref{leAtntoAt}.
Then for any $t>0$ and neighborhood $U$ of $e$ with $\eta(T,\partial U)=0$, as $n\to\infty$,
\begin{eqnarray}
&& \sum_{t_{ni}\leq t}\int_{U^c}F(y^n(t_{n\,i-1}),b_{ni},x)\mu_{ni}(dx) \nonumber \\
&\to& \int_0^t\int_{U^c}F(y(s),e,x)\eta^c(ds,dx) + \sum_{u\leq
t,\,u\in J}\int_{U^c} F(y(u-),h_u,x)\nu_u(dx), \label{fyntni1to}
\end{eqnarray}
and along the subsequence of $n\to\infty$ in Lemma~\ref{leAtntoAt},
\begin{eqnarray}
&& \sum_{t_{ni}\leq t}\sum_{j,k=1}^d\int_G F_{jk}(y^n(t_{n\,i-1}),b_{ni},x)
[\phi_j(x)-\phi_j(b_{ni})][\phi_k(x)-\phi_k(b_{ni})]\mu_{ni}(dx)
\nonumber \\
&\to& \sum_{j,k=1}^d\{\int_0^t F_{jk}(y(s),e,e) dA_{jk}(s) +
\int_0^t\int_GF_{jk}(y(s),e,x)\phi_j(x)\phi_k(x)\eta^c(ds,dx) \nonumber \\
&&\ \ + \sum_{u\leq t,\,u\in J}\int_G
F_{jk}(y(u-),h_u,x)[\phi_j(x)-\phi_j(h_u)]
[\phi_k(x)-\phi_k(h_u)]\nu_u(dx)\}. \label{fjkyntni1to}
\end{eqnarray}
\end{lemm}

\noindent {\bf Proof} \ Let $V$ be a relatively compact neighborhood
of $e$. By Proposition~\ref{pretanbdd}, $\eta_n(t,V^c)\downarrow 0$
as $V\uparrow G$ uniformly in $n$. Because of the
uniform convergence $r(y^n(t_{ni}),y(t_{ni}))\to 0$, $F(y^n(t_{ni}),b,x)-
F(y(t_{ni}),b,x)\to 0$ as $n\to\infty$ uniformly for $t_{ni}\leq t$ and
for $(b,x)$ in a compact set. Therefore, it suffices to prove (\ref{fyntni1to}) with
$y^n$ and $U^c$ replaced by $y$ and $U^c\cap V$ for an arbitrary
relatively compact neighborhood $V$ of $e$ with $\eta(T,\partial
V)=0$. Similarly, it suffices to prove (\ref{fjkyntni1to}) with
$y^n$ and $G$ replaced by $y$ and $V$.

We now show that for $s<t$ and $f\in C_b(G)$ vanishing in a
neighborhood of $e$, as $n\to\infty$,
\begin{equation}
\sum_{s<t_{n\,i-1}<t_{ni}<t}\mu_{ni}(f)\ \to\ \eta((s,\,t)\times f).
\label{sumst}
\end{equation}
The sum $\sum_{s<t_{n\,i-1}<t_{ni}<t} \mu_{ni}(f)$ differs from $\eta_n(t,f)-\eta_n(s,f)$ by one or two
terms of the form $\mu_{ni}(f)$. By Proposition~\ref{pretantoeta}, $\eta_n(t,f)-\eta_n(s,f)\ \to
\ \eta(t,f)-\eta(s,f) = \eta((s,\,t]\times f) = \eta((s,\,t)\times f) +
\nu_t(f)$. Using $J\subset\Delta_n$ as $n\to\infty$,
(\ref{sumst}) may be derived by noting that as $n\to\infty$, $\mu_{ni}(f)\to \nu_t(f)$ if $t\in J\cap
(t_{n\,i-1},\,t_{ni}]$ and $\mu_{ni}(f)\to 0$ otherwise.

Let $\Delta_m'$: $0=s_{m0}<s_{m1}<\cdots <s_{mk}\uparrow\infty$ be a sequence of partitions of $\bR_+$ such
that $\Vert\Delta_m'\Vert\to 0$ as $m\to\infty$. Let $J'$ be the union of $J$ and the set of discontinuity times of $y(t)$.
We may assume $J'\subset\Delta_m'$ as $m\to\infty$, that is, any $v\in J'$ is contained in $\Delta_m'$ for large $m$.
In the following computation,
we will write $A(m,n)\sim B(m,n)$ for any two real valued expressions $A(m,n)$ and $B(m,n)$ if $\vert A(m,n)-B(m,n)\vert\to 0$ as $m\to\infty$ uniformly in $n$.
\begin{eqnarray*}
&& \sum_{s_{mj}\leq t}\{\sum_{s_{m\,j-1}<t_{n\,i-1}<t_{ni}<s_{mj}}
\int_{U^c\cap V}F(y(t_{n\,i-1}),b_{ni},x)\mu_{ni}(dx)\} \\
&\sim& \sum_{s_{mj}\leq
t}\{\sum_{s_{m\,j-1}<t_{n\,i-1}<t_{ni}<s_{mj}} \int_{U^c\cap
V}F(y(s_{m\,j-1}),e,x)\mu_{ni}(dx)\}\ \ \mbox{(because $J'\subset\Delta_m'$ as $m\to\infty$
% ,} \\
% && \mbox{$\forall\varepsilon>0$, first choose $m_0$ so that $\forall m>m_0$, $\vert
% F(y(s_{m\,j-1}),e,x)-F(y(u),x_{v-},x)\vert<\varepsilon$ for} \\
% && \mbox{$u,v\in(s_{m\,j-1},\,s_{mj})$ and $x\in U^c\cap V$, next choose $n_0$ so that $\forall n>n_0$, $\vert
% F(y(s_{m\,j-1}),e,x)-$} \\
% && \mbox{$F(y(t_{n\,i-1}),b_{ni},x)\vert<2\varepsilon$ for $[t_{n\,i-1},\,t_{ni}]\subset (s_{m\,j-1},\,s_{mj})$
% and $x\in U^c\cap V$, then increase $m_0$} \\
% && \mbox{so that $\forall m>m_0$, each $(s_{m\,j-1},\,s_{mj})$ contains
% at most one point of $\Delta_n$ for $n\leq n_0$
)} \\
&\to& \sum_{s_{mj}\leq t}\int_{U^c\cap V}
F(y(s_{m\,j-1}),e,x)\eta((s_{m\,j-1},\,s_{mj}),dx)\ \ \ \
\mbox{(as $n\to\infty$, by (\ref{sumst}))} \\
&\sim& \sum_{s_{mj}\leq t}\int_{(s_{m\,j-1},\,s_{mj})}\int_{U^c\cap V}F(y(s),e,x)\eta^c(ds,dx).
\end{eqnarray*}
On the other hand, for $t_{n\,i-1}<s_{mj}\leq t_{ni}$, Because $\eta(T,\partial U)=\eta(T,\partial V)=0$,
\[\int_{U^c\cap V}F(y(t_{n\,i-1}),b_{ni},x)\mu_{ni}(dx)\ \to\ \nu_{s_{mj}}
(F(y(s_{mj}-),h_{s_{mj}},\cdot)1_{U^c\cap V})\ \ \ \ {\rm as}\
n\to\infty.\] This proves (\ref{fyntni1to}) when $y^n$ and $U^c$ are
replaced by $y$ and $U^c\cap V$. As noted earlier, this proves
(\ref{fyntni1to}) in its original form. The convergence in
(\ref{fjkyntni1to}) is proved in a similar fashion. \ $\Box$
\vspace{2ex}

Now assume $x_0=e$. For $0\leq t<t_{n1}$, let $x_0^n=b_0^n=z_t^0=e$, and with $i\geq 1$, let
\begin{eqnarray}
x_t^n &=& x_{t_{ni}} = x_{n1}x_{n2}\cdots x_{ni},\ \ \ \ t_{ni}\leq t < t_{n\,i+1},
\ \ \ \ \mbox{($x_{ni}=x_{t_{n\,i-1}}^{-1}x_{t_{ni}}$)} \label{xtn} \\
b_t^n &=& b_{n1}b_{n2}\cdots b_{ni},\ \ \ \ t_{ni}\leq t < t_{n\,i+1}, \label{btn} \\
z_t^n &=& x_t^n(b_t^n)^{-1} = z_{n1}z_{n2} \cdots z_{ni},
\ \ \ \ {\rm where}\ \ z_{ni}=b_{t_{n\,i-1}}^n(x_{ni}b_{ni}^{-1})(b_{t_{n\,i-1}}^n)^{-1}. \label{ztn}
\end{eqnarray}

For $f\in C_c^\infty(G)$, let $M_t^nf=f(z_t^n)=f(e)$ for $0\leq
t<t_{n1}$, and for $t\geq t_{n1}$, let
\begin{equation}
M_t^nf = f(z_t^n) - \sum_{t_{ni}\leq t} \int_G [f(z_{t_{n\,i-1}}^n
b_{t_{n\,i-1}}^nxb_{ni}^{-1}(b_{t_{n\,i-1}}^n)^{-1}) -
f(z_{t_{n\,i-1}}^n)]\mu_{ni}(dx). \label{Mtnf}
\end{equation}

\begin{lemm} \label{leMtnmart}
$M_t^nf$ is a martingale under the natural filtration $\myF_t^x$ of process $x_t$.
\end{lemm}

\noindent {\bf Proof} \ The martingale property can be verified directly noting
\begin{eqnarray*}
E[\int_G f(z_{t_{n\,i-1}}^n b_{t_{n\,i-1}}^nxb_{ni}^{-1}
(b_{t_{n\,i-1}}^n)^{-1})\mu_{ni}(dx)\mid\myF_{t_{n\,i-1}}^x] &=&
E[f(z_{t_{n\,i-1}}^n z_{ni})\mid\myF_{t_{n\,i-1}}^x] \\
&=& E[f(z_{t_{ni}}^n)\mid\myF_{t_{n\,i-1}}^x].\ \ \Box
\end{eqnarray*}

Fix a left invariant Riemannian metric $r$ on $G$. A subset of $G$
is relatively compact if and only if it is bounded in $r$. Consider
the following equi-continuity type condition on $b_t^n$.
\vspace{2ex}

\noindent (B) \ For any $T>0$ and $\varepsilon>0$, there are
integers $n_0,m>0$ and $\delta>0$ such that if $n\geq n_0$ and $s,t\in [0,\,T]$
with $0<t-s<\delta$ and $(s,\,t]\cap J_m=\emptyset$, then $r(b_s^n,b_t^n)<\varepsilon$.

\begin{lemm} \label{lebtntobt}
Under (B), $b_t^n$ is bounded for $0\leq t\leq T$ and $n\geq 1$. Moreover, there is a rcll
path $b_t$ in $G$ with $b_0=e$ such that along a subsequence of $n\to\infty$, $b_t^n\to b_t$
uniformly for $0\leq t\leq T$. Furthermore, $b_{t-}^{-1}b_t=h_t$ (the mean of $x_{t-}^{-1}x_t$)
for any $t>0$.
\end{lemm}

\noindent {\bf Proof} \ The boundedness of $b_t^n$ may be derived from (B) in the same way as
the boundedness of $\eta_n(t,U^c)$ in Proposition~\ref{pretanbdd} is derived
from Lemma~\ref{leetanequicont}.

Because $b_t^n$ is bounded, there is a subsequence of $n\to\infty$
along which $b_t^n$ converges to some $b_t$ for $t$ in a countable dense subset
of $[0,\,T]$ including $J$. The equi-continuity of $b_t^n$ in (B) implies
that the convergence holds for all $t$ and is uniform for $t\leq T$, and the limit $b_t$ is
a rcll path with jump $b_{t-}^{-1}b_t=h_t$ at $t\in J$, noting $h_t=\lim_{n\to\infty}b_{ni}$
for $t\in (t_{n\,i-1},\,t_{ni}]$. \ $\Box$

\begin{lemm} \label{leztmart}
Assume (A) and (B), and let $z_t=x_tb_t^{-1}$. Then for $f\in C_c^\infty(G)$, $M_tf$
given in (\ref{Mtf}) is a martingale under the natural filtration $\myF_t^x$ of $x_t$.
\end{lemm}

\noindent {\bf Proof} \ We need only to show that for any $t>0$, $M_t^nf$ is bounded in $n$
and $M_t^nf\to M_tf$.

We will assume, for the time being, that for any $u>0$, $\nu_u$ is small
in the meaning defined before (\ref{phijb}). Then $\mu_{ni}$ are small for $t_{ni}\leq t$ when $n$ is large.

Let $U$ be a neighborhood of $e$ such that $x=e^{\sum_j\phi_j(x)\xi_j}$ for $x\in U$ and $\eta(T,\partial
U)=0$. Let $\sum_{t_{ni}\leq t}$ be the sum in (\ref{Mtnf})
that defines $M_t^nf$ as $f(z_t^n)-\sum_{t_{ni}\leq t}$. A typical term
in $\sum_{t_{ni}\leq t}$ may be written as $\int_G[f(zbxb'^{-1}b^{-1})-f(z)]\mu(dx)$,
where $z=z_{t_{n\,i-1}}^n$, $b=b_{t_{n\,i-1}}^n$, $b'=b_{ni}$ and $\mu=\mu_{ni}$.
By Taylor expansion of $f(zbxb'^{-1}b^{-1})=f(zb\exp(\sum_{j=1}^d\phi_j(x)\xi_j)b'^{-1}b^{-1})$
at $x=b'$,
\begin{eqnarray*}
&& \int_G [f(zbxb'^{-1}b^{-1}) - f(z)]\mu(dx) = \int_{U^c}[\ldots]\mu(dx) +
\int_U[\ldots]\mu(dx) \\
&=& \int_{U^c}[f(zbxb'^{-1}b^{-1}) - f(z)]\mu(dx) + \int_U\{\sum_j f_j(z,b,b')
[\phi_j(x)-\phi_j(b')]\}\mu(dx) \\
&&\ \ + \int_U
\{\frac{1}{2}\sum_{j,k} f_{jk}(z,b,b')[\phi_j(x)-\phi_j(b')][\phi_k(x)-\phi_k(b')]\}
\mu(dx) + \lambda,
\end{eqnarray*}
where
\begin{equation}
f_j(z,b,b') = \frac{\partial}{\partial\phi_j}f(zbe^{\sum_p\phi_p(x)\xi_p}b'^{-1}b^{-1})
\mid_{x=b'}, \label{fjzbbp}
\end{equation}
\begin{equation}
f_{jk}(z,b,b') = \frac{\partial^2}{\partial\phi_j\partial\phi_k}f(zbe^{\sum_p\phi_p(x)\xi_p}b'^{-1}
b^{-1})\mid_{x=b'}, \label{fjkzbbp}
\end{equation}
and the remainder $\lambda$ satisfies $\vert\lambda\vert\leq c_U\mu(\Vert\phi_\cdot-\phi_\cdot(b')\Vert^21_U)$
with constant $c_U\to 0$ as $U\downarrow\{e\}$. Because $b'$ is the mean of the small $\mu=\mu_{ni}$, $\phi_j(b')=\mu(\phi_j)$ and $\int_U
[\phi_j(x)-\phi_j(b')]\mu(dx)=\int_G-\int_{U^c}=0-\int_{U^c}
[\phi_j(x)-\phi_j(b')]\mu(dx)$. This implies that the sum $\sum_{t_{ni}\leq t}$ in (\ref{Mtnf}) is equal to
\begin{eqnarray}
&&  \sum_{t_{ni}\leq t} \int_{U^c}\{f(zbxb'^{-1}b^{-1}) -
f(z)-\sum_j f_j(z,b,b')[\phi_j(x)-\phi_j(b')]\}\mu(dx) \nonumber \\
&&\ \ + \sum_{t_{ni}\leq t} \int_U\{\frac{1}{2}\sum_{j,k} f_{jk}(z,b,b')
[\phi_j(x)-\phi_j(b')][\phi_k(x)-\phi_k(b')]\}\mu(dx) + \Lambda, \label{sumnp}
\end{eqnarray}
where the remainder $\Lambda$ satisfies $\vert\Lambda\vert\leq c_U\sum_{t_{ni}\leq t}\mu(\Vert\phi_\cdot-\phi_\cdot(b')
\Vert^21_U)=c_U q^n(t,U)$.

Note that the two sums in (\ref{sumnp}) are bounded in absolute value
by $c\eta_n(t,U^c)$ and $cq^n(t,U)$ respectively
for some constant $c>0$. Therefore, $\sum_{t_{ni}\leq t}$ and hence $M_t^nf$ are bounded in $n$.

By (\ref{xtn}), $x_{t_{ni}}^n=x_{t_{ni}}$. Because $z_t^n=x_t^n(b_t^n)^{-1}$,
by the uniform convergence $b_t^n\to b_t$, $r(z_{t_{ni}}^n,z_{t_{ni}})=r((b_{t_{ni}}^n)^{-1},
(b_{t_{ni}})^{-1})\to 0$ as $n\to\infty$ uniformly for $t_{ni}\leq t$.

The above fact allows us to apply (\ref{fyntni1to}) in Lemma~\ref{lefjkntofjk} to the first sum
in (\ref{sumnp}). Let $Y=G\times G$, $y=(z,b)$, and $F(y,b',x)=
f(zbxb'^{-1}b^{-1})-f(z)-\sum_j f_j(z,b,b')[\phi_j(x)-\phi_j(b')]$. Define $y^n(t)$
and $y(t)$ by setting $y^n(t)=(z_{t_{ni}}^n,
b_{t_{ni}}^n)$ for $t_{ni}\leq t<t_{n\,i+1}$, and $y(t)=(z_t,b_t)$.
Note that $f_j(z,b,e)=[\Ad(b)\xi_j]f(z)$. It follows that as $n\to\infty$,
the first sum in (\ref{sumnp}) converges to
\begin{eqnarray}
&& \int_0^t\int_{U^c}\{f(z_sb_sxb_s^{-1})-f(z_s)-\sum_j\phi_j(x)[\Ad(b_s)\xi_j]
f(z_s)\}\eta^c(ds,dx) \nonumber \\
\hspace{-1in}+ && \hspace{-0.3in} \sum_{u\leq t,u\in J}\int_{U^c}
\{f(z_{u-}b_{u-}xh_u^{-1}b_{u-}^{-1}) - f(z_{u-}) - \sum_j
f_j(z_{u-},b_{u-},h_u)[\phi_j(x)-\phi_j(h_u)]\}\nu_u(dx). \label{1stsuminsum}
\end{eqnarray}
  Similarly, apply (\ref{fjkyntni1to}) in Lemma~\ref{lefjkntofjk} to
the second sum in (\ref{sumnp}) shows that it converges to
\begin{eqnarray}
&& \frac{1}{2}\sum_{j,k}\int_0^t[\Ad(b_s)\xi_j][\Ad(b_s)\xi_k]f(z_s)dA_s \nonumber \\
&&\ \ + \frac{1}{2}\sum_{j,k}\int_0^t\int_U[\Ad(b_s)\xi_j][\Ad(b_s)\xi_k]f(z_s)\phi_j(x)\phi_k(x)
\eta^c(ds,dx) \nonumber \\
&&\ \ + \frac{1}{2}\sum_{j,k}\sum_{u\leq t,u\in J}\int_Uf_{jk}(z_{u-},b_{u-},h_u)[\phi_j(x)-\phi_j(h_u)]
[\phi_k(x)-\phi_k(h_u)]\nu_u(dx). \label{2ndsuminsum}
\end{eqnarray}

By a computation similar to the one leading to (\ref{sumnp}), using the Taylor expansion of
$f(z_{u-}b_{u-}xh_u^{-1}b_{u-}^{-1})$ at $x=h_u$, one can show
that $\sum_{u\leq t,u\in
J}\int_G[f(z_{u-}b_{u-}xh_u^{-1}b_{u-}^{-1})- f(z_{u-})]\nu_u(dx)$ is the sum
of $\sum_{u\leq t,u\in J}(\cdots)$ in (\ref{1stsuminsum}) and $(1/2)\sum_{j,k}\sum_{u\leq t,u\in J}(\cdots)$
in (\ref{2ndsuminsum}), plus an error term that converges to $0$ as $U\downarrow\{e\}$.
Letting $U\downarrow\{e\}$ in (\ref{1stsuminsum}) and (\ref{2ndsuminsum}) shows that $M_t^nf\to M_tf$.

For $u\in J$, let $a_u$ be the $\nu_u$-integral term in $M_tf$ given by (\ref{Mtf})
and let $a_u^n$ be the $\mu_{ni}$-integral term in $M_t^nf$ given by (\ref{Mtnf})
with $u\in (t_{n\,i-1},\,t_{ni}]$. Then $a_u^n\to a_u$ as $n\to\infty$.
In the preceding computation, we have assumed that all $\nu_u$ are small.
When $\nu_u$ is not small, $\mu_{ni}$ with $u\in (t_{n\,i-1},\,t_{ni}]$ may not be small
even for large $n$, so the computation leading to (\ref{sumnp}) done for the $\mu_{ni}$-integral is
not valid, but we may still take limit $a_u^n\to a_u$ for this term and the result is
the same. Since there are only finitely many $u\in J$ for which $\nu_u$ are not small, the
result holds even when not all $\nu_u$ are small. \ $\Box$

\section{Proofs of (A) and (B)} \label{secAB}

Because $q^n(t,\cdot)=q^n(t_{ni},\cdot)$ and $b_t^n=b_{t_{ni}}^n$ if $t_{ni}\leq
t<t_{n\,i+1}$, and $\Delta_n$ contains $J$ as $n\to\infty$, to verify (A) and (B),
we may assume $s,t$ in (A) and (B) are contained in $[0,\,T]\cap\Delta_n$.

If either (A) or (B) does not hold, then for some $\varepsilon>0$ and any $n_0,m,\delta>0$, there are $s_n,t_n\in[0,\,T]\cap\Delta_n$
with $n\geq n_0$, $0<t_n-s_n<\delta$ and $(s_n,\,t_n]\cap J_m=\emptyset$ such
that either $q^n(t_n,G)-q^n(s_n,G)\geq\varepsilon$ or $r(b_{s_n}^n,
b_{t_n}^n)\geq\varepsilon$. Because the jumps of $q^n(t,G)$ and $b_t^n$ for $t\in (s_n,\,t_n]$ become
arbitrarily small when $n$ and $m$ are large, by decreasing $t_n$ if necessary,
we may also assume $q^n(t_n,G)-q^n(s_n,G)\leq 2\varepsilon$ and $r(b_{s_n}^n,b_t^n)\leq 2\varepsilon$ for $s_n\leq t\leq t_n$.

Letting $\delta\to 0$ and $m\to\infty$ yields a subsequence of $n\to\infty$ such
that $s_n$ and $t_n$ in $[0,\,T]\cap\Delta_n$ converge to a common limit, $(s_n,\,t_n]\cap
J_{m_n}=\emptyset$ with $m_n\uparrow\infty$, and either
\vspace{2ex}

\noindent (i) \ $\varepsilon\leq q^n(t_n,G)-q^n(s_n,G)\leq
2\varepsilon$ and $r(b_{s_n}^n,b_t^n)\leq 2\varepsilon$ for $t\in
[s_n,\,t_n]$, or

\noindent (ii) \ $q^n(t_n,G)-q^n(s_n,G)\leq 2\varepsilon$, $r(b_{s_n}^n,b_{t_n}^n)
\geq\varepsilon$ and $r(b_{s_n}^n,b_t^n)\leq
2\varepsilon$ for $t\in [s_n,\,t_n]$.
\vspace{2ex}

We will derive a contradiction from (i) or (ii). Let $\varepsilon_n=\sup_{s_n<t\leq t_n}
[q^n(t,G)-q^n(t-,G)]$. Then $\varepsilon_n\to 0$ as $n\to\infty$. Let $\gamma_n$: $[0,\,\infty)
\to[s_n,\,\infty)$ be a strictly increasing continuous function such that $\gamma_n([0,\,1])= [s_n,\,t_n]$.
Because $q^n(t,G)-q^n(s_n,G)$ is a nondecreasing step function in $t$, bounded by $2\varepsilon$ on $[s_n,\,t_n]$
with jump size $\leq\varepsilon_n$, $\gamma_n$ may be chosen so
that $\vert q^n(\gamma_n(t),G)-q^n(\gamma_n(s),G)\vert \leq
2\varepsilon\vert t-s\vert+\varepsilon_n$ for $s,t\in [0,\,1]$.

Let $x_t^{\gamma_n}=x_{s_n}^{-1}x_{\gamma_n(t)}$ for $t\leq 1$ and $x_t^{\gamma_n}=x_{s_n}^{-1}x_{t_n}$ for $t>1$.
Then $x_t^{\gamma_n}$ is an inhomogeneous \levy process in $G$ starting at $e$, and on $[0,\,1]$, it is just
the process $x_t$ on $[s_n,\,t_n]$ time changed by $\gamma_n$.
Because $s_n$ and $t_n$ converge to a common limit, and $(s_n,\,t_n]\cap J_{m_n}=\emptyset$
with  $m_n\uparrow\infty$, it follows that $x_t^{\gamma_n}\to e$ as $n\to\infty$ uniformly in $t$ a.s.. Define
\[x_t^{\gamma,n},\ \mu_{ni}^\gamma,\ \eta_n^\gamma(t,\cdot),\ b_{ni}^\gamma,
\ b_t^{\gamma,n},\ A^{n,\gamma}(t,\cdot),\ q^{\gamma,n}(t,\cdot),
\ z_t^{\gamma,n},\ M_t^{\gamma,n}f,\] for the time changed process $x_t^{\gamma_n}$ and
partition $\Delta_n^\gamma=\{s_{ni}\}$, where $s_{ni}=\gamma_n^{-1}(t_{ni})$, in the same way as $x_t^n,\mu_{ni},
\eta_n(t,\cdot),b_{ni},b_t^n$, $A^n(t,\cdot),q^n(t,\cdot),z_t^n, M_t^nf$ are defined for the process $x_t$
and partition $\Delta_n=\{t_{ni}\}$. Then for $s,t\in [0,\,1]$,
\begin{equation}
\vert q^{\gamma,n}(t,G) - q^{\gamma,n}(s,G)\vert = \vert
q^n(\gamma_n(t),G)-q^n(\gamma_n(s),G)\vert\ \leq\ 2\varepsilon
\vert t-s\vert + \varepsilon_n. \label{qgamman}
\end{equation}
Because $\varepsilon_n\to 0$, the above means that $q^{\gamma,n}(t,G)$ are equi-continuous in $t$ for large $n$.

Note that we now have a process $x_t^{\gamma_n}$ for each $n$ and
from which other objects, such as $\eta_n^\gamma$ and $A^{\gamma,n}$, are defined,
unlike before when we have a single process $x_t$ for all $n$, but as $x_t^{\gamma_n}\to e$
uniformly in $t$, the results established for $\eta_n$ and $A^n$ hold for $\eta_n^\gamma$
and $A^{\gamma,n}$ in simpler forms. For example, because $\eta_n^\gamma(t,\cdot)=
\eta_n(\gamma_n(t),\cdot)-\eta_n(s_n,\cdot)$, by Proposition~\ref{pretantoeta}, $\eta_n^\gamma(t,f)\to 0$ as $n\to\infty$
for $f\in C_b(G)$ vanishing in near $e$. Using
(\ref{qgamman}), the proofs of Lemmas \ref{leAtntoAt} and
\ref{lefjkntofjk} can be easily modified for $A^{\gamma,n}$, and also simplified,
to show that there is a covariance matrix
function $A^\gamma(t)$ such that $A^{\gamma,n}(t,f)\to
f(e)A^\gamma(t)$ for any $f\in C_b(G)$, and under the assumption of
Lemma~\ref{lefjkntofjk},
\begin{equation}
\sum_{s_{ni}\leq t}\int_{U^c}F(y^n(s_{n\,i-1}),b_{ni}^\gamma,x)
\mu_{ni}^\gamma(dx)\ \to\ 0 \label{fyntni1togamma}
\end{equation}
and
\begin{eqnarray}
&& \sum_{s_{ni}\leq t}\sum_{j,k=1}^d F_{jk}(y^n(s_{n\,i-1}),
b_{ni}^\gamma,x)[\phi_j(x)-\phi_j(b_{ni}^\gamma)]
[\phi_k(x)-\phi_k(b_{ni}^\gamma)]\mu_{ni}^\gamma(dx) \nonumber \\
&\to& \sum_{j,k=1}^d\int_0^t F_{jk}(y(s),e,e)dA_{jk}^\gamma(s). \label{fjkyntni1togamma}
\end{eqnarray}

Let $D(G)$ be the space of rcll paths $\bR_+\to G$. Equipped with
the Skorohod metric, $D(G)$ is a complete separable metric space (see \cite[chapter~3]{ek}).
A rcll process $x_t$ in $G$ may be regarded as a random variable in $D(G)$. A sequence
of processes $x_t^n$ is said to converge weakly in $D(G)$ to a process $x_t$ if
their distributions converge weakly on $D(G)$ to the distribution of process $x_t$.

We will show that $z_t^{\gamma,n}$ has a weakly convergent subsequence in $D(G)$. Let $U$ be a neighborhood of $e$ and
let $\sigma$ be a stopping time. The amount of time it takes for a process $x_t$ to make a $U^c$-displacement from time $\sigma$ is denoted as $\tau_U^\sigma$,
that is,
\begin{equation}
\tau_U^\sigma = \inf\{t>0;\ \ x_\sigma^{-1}x_{\sigma+t}\in U^c\}
\label{tauUsigma}
\end{equation}
For a sequence of processes $x_t^n$ in $G$, let $\tau_U^{\sigma,n}$
be the $U^c$-displacement time for $x_t^n$ from $\sigma$.

\begin{lemm} \label{lerelweakcompt}
A sequence of processes $x_t^n$ in $G$ have a weakly convergent subsequence
in $D(G)$ if for any $T>0$ and any neighborhood $U$ of $e$,
\begin{equation}
\overline{\lim_{n\to\infty}}\sup_{\sigma\leq T} P(\tau_U^{\sigma,n}
<\delta)\ \to\ 0\ \ {\rm as}\ \delta\to 0, \label{supPtauU}
\end{equation}
and
\begin{equation}
\overline{\lim_{n\to\infty}}\sup_{\sigma\leq T} P[(x_{\sigma-}^n)^{-1} x_\sigma^n\in K^c]\ \to\ 0 \ \ \mbox{as compact}\ K\uparrow G
\ \ \mbox{(in the topology on $G$),} \label{supPxsigma}
\end{equation}
where $\sup_{\sigma\leq T}$ is taken over all stopping times
$\sigma\leq T$.
\end{lemm}

\noindent {\bf Proof} \ This lemma is proved in \cite{feinsilver}.
We will provide a different proof as the argument will be used to prove Lemma~\ref{lerelweakcompt2}.
Let $r$ be a left invariant metric on $G$ as before. As in \cite[section~3.6]{ek}, the measurement of $\delta$-oscillation of a path $x$ in $D(G)$
on $[0,\,T]$ is given by
\begin{equation}
\omega'(x,\delta,T) = \inf_{\{t_i\}}\max_{1\leq i\leq n}\sup_{s,t\in
[t_{i-1},\,t_i)}r(x_t,x_s), \label{omegaprime}
\end{equation}
where the infimum $\inf_{\{t_i\}}$ is taken over all partitions
$0=t_0<t_1<\cdots<t_{n-1}<T\leq t_n$ with $\min_{1\leq i\leq
n}(t_i-t_{i-1})>\delta$. By Corollary~7.4 in \cite[chapter~3]{ek}, $x_t^n$ have
a weakly convergent subsequence in $D(G)$ if for any $T>0$,
\begin{equation}
\overline{\lim_{n\to\infty}}P(x_t^n\in K^c\ {\rm for\ some}\ t\leq
T)\ \to\ 0\ \ {\rm as\ compact}\ K\uparrow G\ \ \mbox{(in the $r$-topology)} \label{xtnKc}
\end{equation}
and for any $\varepsilon>0$
\begin{equation}
\overline{\lim_{n\to\infty}}P[\omega'(x_\cdot^n,\delta,T)\geq\varepsilon]\
\to\ 0 \ \ {\rm as}\ \delta\to 0. \label{omegapvarepsilon}
\end{equation}

For a fixed $\varepsilon>0$, the successive stopping times
$0=\tau_0^\varepsilon<\tau_1^\varepsilon<
\tau_2^\varepsilon<\cdots<\tau_i^\varepsilon<\cdots$ when a rcll
process $x_t$ makes an $\varepsilon$-displacement are defined
inductively by
\[\tau_i^\varepsilon = \inf\{t>\tau_{i-1}^\varepsilon;\ \ r(x_t,
x_{\tau_{i-1}^\varepsilon})>\varepsilon\}\]
for $i=1,2,3,\ldots$, setting $\inf\emptyset=\infty$ and $\tau_i^\varepsilon=\infty$
if $\tau_{i-1}^\varepsilon=\infty$. Let $\tau_i^{\varepsilon,n}$ be
the $\varepsilon$-displacement times of the process $x_t^n$.
It is easy to see that
\begin{equation}
[\min_{i\geq 0}\{\tau_{i+1}^\varepsilon-\tau_i^\varepsilon;
\,\tau_i^\varepsilon<T\}\ >\ \delta]\ \ \ \ {\rm implies} \ \ \ \
[\omega'(x,\delta,T)\leq 2\varepsilon], \label{mintauomega}
\end{equation}
and hence $P[\omega'(x,\delta,T)
>2\varepsilon]\leq P[\min_{i\geq 0}\{\tau_{i+1}^\varepsilon-
\tau_i^\varepsilon;\,\tau_i^\varepsilon<T\}\leq\delta]$. Thus, if
\begin{equation}
\forall T>0\ {\rm and}\ \varepsilon>0,\ \ \ \ \overline{
\lim_{n\to\infty}} P[\min_{i\geq 0}\{\tau_{i+1}^{\varepsilon,n}-
\tau_i^{\varepsilon,n}; \,\tau_i^{\varepsilon,n}<T\}\ <\ \delta]\
\to\ 0\ \ {\rm as}\ \delta\to 0, \label{supPmintau}
\end{equation}
then (\ref{omegapvarepsilon}) holds for any $T>0$ and $\varepsilon>0$.

Let $F(t)=\sup_n\sup_{i\geq 0}P(\tau_{i+1}^{\varepsilon,n}
-\tau_i^{\varepsilon,n}<t,\,\tau_i^{\varepsilon,n}<T)$. By Lemma~8.2
in \cite[chapter~3]{ek},
\begin{equation}
F(\delta)\ \leq\ \sup_n P[\min_{i\geq
0}\{\tau_{i+1}^{\varepsilon,n}- \tau_i^{\varepsilon,n};
\,\tau_i^{\varepsilon,n}<T\}\ <\ \delta]\ \leq\ LF(\delta) +
\int_0^\infty e^{-Lt}F(t/L)dt \label{Fdelta}
\end{equation}
for any $\delta>0$ and $L=1,2,3,\ldots$. It is clear that $\sup_n$
in $F(t)$ and (\ref{Fdelta}) may be replaced by $\lims_{n\to\infty}$. Consequently,
(\ref{supPmintau}) is equivalent to
\begin{equation}
\forall T>0\ {\rm and}\ \varepsilon>0,\ \ \ \ \overline{
\lim_{n\to\infty}}\sup_{i\geq 0}
P(\tau_{i+1}^{\varepsilon,n}-\tau_i^{\varepsilon,n}<\delta;\
\tau_i^{\varepsilon,n}<T)\ \to\ 0\ \ {\rm as}\ \delta\to 0.
\label{supsupPtau}
\end{equation}

It is now clear that (\ref{supsupPtau}), and hence
(\ref{supPmintau}) and (\ref{omegapvarepsilon}), are implied by
(\ref{supPtauU}) for any neighborhood $U$ of $e$ and $T>0$. It
remains to verify (\ref{xtnKc}). It suffices to show that for any
$\eta>0$, there are a compact $K\subset G$ and integer $m>0$ such
that $P[x_t^n\in K^c$ for some $t\leq T]\leq\eta$
for all $n\geq m$. Because (\ref{supPtauU}) implies
(\ref{supPmintau}), there are $\delta>0$ and integer $m>0$ such that
for all $n\geq m$, $P(A_n)<\eta/2$, where $A_n=[\min_{i\geq
0}\{\tau_{i+1}^{\varepsilon,n}-
\tau_i^{\varepsilon,n};\,\tau_i^{\varepsilon,n}<T\}<\delta]$. Let $p=[T/\delta]$,
the integer part of $T/\delta$. Then on $A_n^c$, $\min_{i\geq 0}\{\tau_{i+1}^{\varepsilon,n}-
\tau_i^{\varepsilon,n};\,\tau_i^{\varepsilon,n}<T\} \geq\delta$, and
hence $x_t^n$ makes at most $p$ displacements of size $\varepsilon$,
before time $T$, between $\tau_1^{\varepsilon,n},\tau_2^{\varepsilon,n},\ldots, \tau_p^{\varepsilon,n}$,
with possible jumps at these times. By (\ref{supPxsigma}),
there is a compact $H\subset G$ such that $P(B_{n,i})\leq\eta/(2p)$
for all $n\geq m$ and $1\leq i\leq p$, where
$B_{n,i}=[(x_{\tau-}^n)^{-1} x_\tau^n\in H^c;\,\tau\leq T]$ with
$\tau= \tau_i^{\varepsilon,n}$. Let $U$ be the $\varepsilon$-ball
around $e$ and let $K$ be a compact subset of $G$ containing
$UHUH\cdots UHU= \{u_1h_1u_2h_2\cdots u_ph_pu_{p+1}$; $u_i\in U$ and
$h_i\in H\}$. Then for $n\geq m$,
\[P[x_t^n\in K^c\ \mbox{for some $t\leq T$}]\ \leq
\ P(A_n) + \sum_{i=1}^pP(A_n^c\cap B_{n,i})\ \leq\ \eta. \ \ \Box\]

We will now apply Lemma~\ref{lerelweakcompt} to the processes $z_t^{\gamma,n}$.
Because $z_t^{\gamma,n}$ is constant for $t\geq
1$, to verify the conditions of Lemma~\ref{lerelweakcompt}, it is
enough to consider only stopping times $\sigma\leq 1$.
Let $f\in C_c^\infty(G)$ be such that $0\leq f\leq 1$ on $G$, $f(e)=1$ and $f=0$ on $U^c$.
Write $\tau$ for the $U^c$-displacement time for process $z_t^{\gamma,n}$
from a stopping time $\sigma$. Let $f_\sigma=f\circ (z_\sigma^{\gamma,n})^{-1}$. Then
\[P(\tau<\delta) = E[f_\sigma(z_\sigma^{\gamma,n})-
f_\sigma(z_{\sigma+\tau}^{\gamma,n});\,\tau<\delta] \leq
E[f_\sigma(z_\sigma^{\gamma,n})-
f_\sigma(z_{\sigma+\tau\wedge\delta}^{\gamma,n})],\]
noting $f_\sigma(z_\sigma^{\gamma,n})=1$, $f_\sigma(z_{\sigma+\tau}^{\gamma,n})=0$
and $\tau=\tau\wedge\delta$ on $[\tau<\delta]$, where $a\wedge b=\min(a,b)$. Because
\begin{equation}
M_t^{\gamma,n}f = f(z_t^{\gamma,n})-\sum_{s_{ni}\leq t}
\int_G[f(z_{s_{n\,i-1}}^{\gamma,n}b_{s_{n\,i-1}}^{\gamma,n}x(b_{ni}^{\gamma})^{-1}
(b_{s_{n\,i-1}}^{\gamma,n})^{-1})-
f(z_{s_{n\,i-1}}^{\gamma,n})]\mu_{ni}^\gamma(dx) \label{Mtgammanf}
\end{equation}
is a martingale for any $f\in C_c^\infty(G)$, and $\sigma$ and
$\sigma+\tau\wedge\delta$ are stopping times,
\[E[M_\sigma^{\gamma,n}f_\sigma - M_{\sigma+\tau\wedge\delta}^{\gamma,n}
f_\sigma] = E\{E[M_\sigma^{\gamma,n}f_\sigma -
M_{\sigma+\tau\wedge\delta}^{\gamma,n} f_\sigma\mid\myF_\sigma]\} =
0.\]
Writing $z,b,b',\mu$ for $z_{s_{n\,i-1}}^{\gamma,n},
b_{s_{n\,i-1}}^{\gamma,n},b_{ni}^\gamma, \mu_{ni}^\gamma$,
we obtain
\begin{eqnarray}
P(\tau<\delta) &\leq& - E\{\sum_{\sigma<s_{ni}\leq
\sigma+\tau\wedge\delta}\int_G
[f_\sigma(zbxb'^{-1}b^{-1})-f_\sigma(z)]\mu(dx)\} \nonumber \\
&\leq& E\{\sum_{\sigma<s_{ni}\leq\sigma+\delta}\vert
\int_G[f_\sigma(zbxb'^{-1}b^{-1})-f_\sigma(z)]\mu(dx)\vert\}.
\label{Ptaudelta}
\end{eqnarray}
Performing the same computation with Taylor's expansion as in the
proof of Lemma~\ref{leztmart} and noting that $\mu$ is small for large $n$,
\begin{eqnarray}
&& \int_G[f(zbxb'^{-1}b^{-1})-f(z)]\mu(dx) \nonumber \\
&=& \int_{U^c}[f(zbxb'^{-1}b^{-1}) -
f(z)-\sum_j f_j(z,b,b')(\phi_j(x)-\phi_j(b'))]\mu(dx) \nonumber \\
&&\ \ + \int_U[\frac{1}{2}\sum_{j,k} f_{jk}(z,b,b')
(\phi_j(x)-\phi_j(b'))(\phi_k(x)-\phi_k(b'))]\mu(dx) + \lambda,
\label{intfzbxmu}
\end{eqnarray}
where the remainder $\lambda$ satisfies $\vert\lambda\vert\leq c_U\mu
(\Vert\phi_\cdot- \phi_\cdot(b')\Vert^21_U)$ with constant $c_U\to 0$ as $U\downarrow\{e\}$.
By (\ref{Ptaudelta}) and (\ref{intfzbxmu}), $P(\tau<\delta)$ is controlled by
\[E[\eta_n^\gamma(1+\delta,U^c) + q^{\gamma,n}(\sigma+\delta,G) -
q^{\gamma,n}(\sigma,G)].\] Because $\eta_n^\gamma(t,\cdot)\to 0$
weakly outside a neighborhood of $e$ and $q^{\gamma,n}(t,G)$ is
equi-continuous in $t$ for large $n$ as expressed by
(\ref{qgamman}), it follows that $\lims_{n\to\infty}\sup_{\sigma\leq
1}P(\tau<\delta)\to 0$ as $\delta\to 0$. This verifies the condition (\ref{supPtauU}) in
Lemma~\ref{lerelweakcompt} for $x_t^n=z_t^{\gamma,n}$.

To verify (\ref{supPxsigma}), note that because $x_t^{\gamma,n}=z_t^{\gamma,n}
b_t^{\gamma,n}$,
\[P[(z_{\sigma-}^{\gamma,n})^{-1}z_\sigma^{\gamma,n}\in K^c] =
P[(x_{\sigma-}^{\gamma,n})^{-1}x_\sigma^{\gamma,n} \in
((b_{\sigma-}^{\gamma,n})^{-1}K b_\sigma^{\gamma,n})^c].\]
Because by either (i) or (ii), $b_t^{\gamma,n}$ are bounded in $n$, when $K$ is
large, $(b_{\sigma-}^{\gamma,n})^{-1}K b_\sigma^{\gamma,n}$ contains a fixed
neighborhood $H$ of $e$. It follows that
\[P[(z_{\sigma-}^{\gamma,n})^{-1}z_\sigma^{\gamma,n}\in K^c] \leq
P[(x_{\sigma-}^{\gamma,n})^{-1}x_\sigma^{\gamma,n} \in H^c] \leq
\eta_n^\gamma(1,H^c) \to 0\] as $n\to\infty$. This establishes
(\ref{supPxsigma}) even before taking $K\uparrow G$.

Now by Lemma~\ref{lerelweakcompt}, along a subsequence of $n\to\infty$, $z_t^{\gamma,n}$ converges
weakly to a rcll process $z_t^\gamma$ in $D(G)$. Because $x_t^{\gamma,n}\to e$ uniformly in $t$, it follows
that $b_\cdot^{\gamma,n}=(z_\cdot^{\gamma,n})^{-1}x_\cdot^{\gamma,n}$ converges to $b_\cdot^\gamma=
(z_\cdot^\gamma)^{-1}$ in $D(G)$. In particular, this means that $z_t^\gamma=(b_t^\gamma)^{-1}$ is non-random.

The convergence $b_t^{\gamma,n}\to b_t^\gamma$ under the Skorohod metric means (see
Proposition~5.3(c) in \cite[chapter~3]{ek}) that there
are continuous strictly increasing functions $\lambda_n$: $\bR_+\to\bR_+$ such
that as $n\to\infty$, $\lambda_n(t)-t\to 0$ and $r(b_t^{\gamma,n},
b_{\lambda_n(t)}^\gamma)\to 0$ uniformly for $0\leq t\leq 1$. If $b_t^\gamma$ has a jump
of size $r(b_{s-}^\gamma,b_s)>0$ at time $s$, then $b_t^{\gamma,n}$ would
have a jump of size close to $r(b_{s-}^\gamma,b_s)$ at time $t=\lambda_n^{-1}(s)$,
which is impossible because the jumps of $b_t^{\gamma,n}$ are uniformly small
when $n$ is large. It follows that $b_t^\gamma$ is continuous in $t$ and
hence $b_t^{\gamma,n}\to b_t^\gamma$ uniformly in $t$ as $n\to\infty$. Then
the convergence $z_t^{\gamma,n}\to z_t^\gamma$ is also uniform in $t$.

By (\ref{fyntni1togamma}), (\ref{fjkyntni1togamma}) and (\ref{intfzbxmu}),
the martingale $M_t^{\gamma,n}f$ given by (\ref{Mtgammanf}) converges to
\begin{equation}
f(z_t^\gamma) - \int_0^t\frac{1}{2}\sum_{j,k=1}^d[\Ad(b_s^\gamma)\xi_j]
[\Ad(b_s^\gamma)\xi_k]f(z_s^\gamma)dA_{jk}^\gamma(s) \label{fztmart}
\end{equation}
for any $f\in C_c^\infty(G)$. It follows that $e=z_t^\gamma
b_t^\gamma$ and for any $f\in C_c^\infty(G)$, (\ref{fztmart}) is a martingale.
This provides a representation of the trivial \levy
process $x_t=e$ by the \levy triple $(b^\gamma,A^\gamma,0)$.
By the uniqueness of the \levy triple, established in \cite{feinsilver}
(also see Lemma~\ref{leuniquetriple} later), $b_t^\gamma=e$ and $A^\gamma(t)=0$. This is
a contradiction because if (i) holds, then $\trace A^\gamma(1)=
\lim_n q^n(t_n,G)\geq\varepsilon$, and if (ii) holds, then
$b_1^\gamma=\lim_n b_{t_n}^n$ is at lease $\varepsilon$ distance
away from $e$. We now have proved that (A) and (B) must hold.

Because (A) and (B) have been verified, by Proposition~\ref{pretaextlevy} and Lemma~\ref{leztmart},
we have proved that an inhomogeneous \levy process $x_t$ is represented $(b,A,\eta)$ with $\eta(t,\cdot)$
being its jump intensity measure, and $A(t)$ and $b_t$ given in Lemmas \ref{leAtntoAt} and
\ref{lebtntobt}. This is the first part of Theorem~\ref{thmartrep}
except for the uniqueness of $(b,A,\eta)$.

\section{Uniqueness of the triple} \label{secuniquetriple}

Let $(b,A,\eta)$ be a \levy triple on $G$ with $b_t$ of finite
variation, and let $\nu=\{\nu_t$; $t\geq 0\}$ be a family of
probability measures on $G$ such that $\nu_t=\delta_e$ except for countably many $t>0$, and
for any $t>0$ and neighborhood $U$ of $e$,
\begin{equation}
  \sum_{u\leq t}\nu_u(U^c) < \infty,\ \ \ \ \sum_{u\leq t}\sum_{j=1}^d\vert\nu_u(\phi_j)\vert <
\infty\ \ {\rm and} \ \ \sum_{u\leq t}\nu_u(\Vert\phi_\cdot\Vert^2)<\infty. \label{sumnu2}
\end{equation}

A rcll process $z_t$ in $G$ is said to have the martingale property
under the quadruple $(b,A,\eta,\nu)$ as described above, or the
$(b,A,\eta,\nu)$-martingale property, if for $f\in C_c^\infty(G)$,
\begin{eqnarray}
&& f(z_t) - \int_0^t\sum_i\xi_if(z_s)db_i(s) -
\int_0^t\frac{1}{2}\sum_{j,k} \xi_j\xi_kf(z_s)dA_{jk}(s) -
\int_0^t\int_G[f(z_sx)-f(z_s) \nonumber \\
&& -\sum_i\phi_i(x) \xi_if(z_s)]\eta(ds,dx) - \sum_{u\leq t}\int_G
[f(z_{u-}x)-f(z_{u-})]\nu_u(dx) \label{ztmartprop}
\end{eqnarray}
is a martingale under the natural filtration $\myF_t^z$ of process $z_t$,
where $b_i(t)$ are the components of $b_t$. In the sequel, a $(b,A,\eta,\nu)$-martingale
property always refers to a \levy triple $(b,A,\eta)$ with $b_t$ of finite variation and a family $\nu$
of probability measures $\nu_t$ on $G$ satisfying (\ref{sumnu2}).

Recall that $[\Ad(g)]$ is the matrix representing $\Ad(g)$ under the
basis $\{\xi_1,\ldots,\xi_d\}$ of $\myg$. Let $[\Ad(g)]'$ be its transpose.

\begin{lemm} \label{leztmartproperty}
If $x_t$ is an inhomogeneous \levy process in $G$ represented by an
extended \levy triple $(b,A,\eta)$, and $x_t=z_tb_t$, then $z_t$ has
the $(\bar{b},\bar{A},\bar{\eta},\bar{\nu})$-martingale property,
where $d\bar{A}(t)=[\Ad(b_t)]dA(t)[\Ad(b_t)]'$, $\bar{b}_t$ has components
\begin{equation}
\bar{b}_i(t) = \int_G\{\phi_i(b_txb_t^{-1})-\sum_p
\phi_p(x)[\Ad(b_t)]_{ip}\}\eta^c(t,dx), \label{barbi}
\end{equation}
$\bar{\eta}(t,f)=\int_G f(b_txb_t^{-1})\eta^c(t,dx)$ and $\bar{\nu}_t(f)= \int_G
f(b_{t-}xh_t^{-1}b_{t-}^{-1})\nu_t(dx)$, where $\eta^c$ is the continuous part of $\eta$
and $\nu_t$ is given in (\ref{nut}) with mean $h_t$.
Thus, $\bar{\nu}_t=\delta_e$ if $\eta(t,\cdot)$ is continuous at time $t$.
\end{lemm}

\noindent {\bf Proof} \ Formally this follows directly from (\ref{Mtf}), but
we need to verify that $\bar{b}_i(t)$ have finite variation and (\ref{sumnu2})
holds for $\bar{\nu}_t$. Because $e^{\sum_j\phi_j(b_uxb_u^{-1})\xi_j}=
b_uxb_u^{-1}=e^{\sum_j\phi_j(x) \Ad(b_u)\xi_j}$ for $u\leq t$
and $x$ in a neighborhood $U$ of $e$, the integrand in (\ref{barbi})
vanishes in $U$. Noting $\eta^c(t,U^c)<\infty$, it is now easy to show that $\bar{b}_i(t)$
have finite variation. The first inequality in (\ref{sumnu2}) for $\bar{\nu}$, $\sum_{u\leq t}
\bar{\nu}_u(U^c)<\infty$, follows from $\sum_{u\leq t}\nu_u(U^c)<\infty$.
The second and the third inequalities in (\ref{sumnu2}), $\sum_{u\leq t}\vert\bar{\nu}_u(\phi_j)
\vert=\sum_{u\leq t}\vert\int\phi_j(b_{u-}xh_u^{-1}b_{u-}^{-1})\nu_u(dx)\vert<\infty$
and $\sum_{u\leq t}\bar{\nu}_u(\Vert\phi_\cdot\Vert^2)=\sum_{u\leq t}\int\Vert\phi_\cdot
(b_{u-}xh_u^{-1}b_{u-}^{-1})\Vert^2\nu_u(dx)<\infty$, follow from a computation similar
to the proof of Lemma~\ref{leztmart} using Taylor expansion of $\phi_j(b_{u-}xh_u^{-1}b_{u-}^{-1})=
\phi_j(b_{u-}e^{\sum_p\phi_p(x)\xi_p}h_u^{-1}b_{u-}^{-1})$ at $x=h_u$. \ $\Box$
\vspace{2ex}

The following lemma summarizes the martingale characterization and
uniqueness theorems in \cite[part~3]{feinsilver} (pages 80-81).

\begin{lemm} \label{lefein}
Given a \levy triple $(b,A,\eta)$ with $b_t$ of finite variation,
there is an inhomogeneous \levy process $z_t$ in $G$ with $z_0=e$, unique
in distribution, such that it has the $(b,A,\eta,0)$-martingale
property. Moreover, given a rcll process $z_t$ in $G$,
there is at most one \levy triple $(b,A,\eta)$ with $b_t$ of finite variation such
that the $(b,A,\eta,0)$-martingale property holds.
\end{lemm}

\begin{lemm} \label{leztfixedjumpeta}
If $z_t$ is a rcll process in $G$ having a $(b,A,\eta,\nu)$-martingale property,
then the fixed jumps of $z_t$ are determined by $\nu$, that is,
the distribution of $z_{t-}^{-1}z_t$ is $\nu_t$ for all $t>0$.
\end{lemm}

\noindent {\bf Proof} \ Let $M_t$ be the martingale
(\ref{ztmartprop}). Then $E[M_t-M_{t-}\mid\myF_{t-}^z]=0$. By
(\ref{ztmartprop}),
\[\forall f\in C_c^\infty(G),\ \ \ \ E\{f(z_t)-f(z_{t-})-\int_G
[f(z_{t-}x)-f(z_{t-})]\nu_t(dx)\mid\myF_{t-}^z\}=0.\]
Then $E[f(z_t)\mid\myF_{t-}^z]=\int_Gf(z_{t-}x)\nu_t(dx)$, and
hence $\nu_t$ is the distribution of $z_{t-}^{-1}z_t$. \ $\Box$

\begin{coro} \label{coxtfixedjumpeta}
If $x_t$ is a rcll process in $G$ represented by an extended \levy
triple $(b,A,\eta)$, then the fixed jumps of $x_t$ are determined by
the discontinuous part of $\eta$, that is, $x_{t-}^{-1}x_t$ has distribution $\nu_t=
\eta(\{t\}\times\cdot)+[1-\eta(\{t\}\times G)]\delta_e$ for any $t>0$.
\end{coro}

\noindent {\bf Proof} \ Let $x_t=z_tb_t$. By
Lemma~\ref{leztmartproperty}, $z_t$ has the $(\bar{b},\bar{A},\bar{\eta},\bar{\nu})$-martingale property, and by
Lemma~\ref{leztfixedjumpeta}, the fixed jumps of $z_t$ are determined by $\bar{\nu}$. Then
\[E[f(x_{t-}^{-1}x_t)] = E[f(b_{t-}^{-1}z_{t-}^{-1}z_tb_t)] =
\int_G f(b_{t-}^{-1}xb_t)\bar{\nu}_t(dx) = \int_G
f(xh_t^{-1}b_{t-}^{-1}b_t)\nu_t(dx) = \nu_t(f)\] for $f\in\myB_+(G)$,
because $b_{t-}^{-1}b_t=h_t$. \ $\Box$
\vspace{2ex}

A possible jump of a rcll process $z_t$ at time $u$ may be removed to obtain a new
process $z_t'$ defined by $z_t'=z_t$ for $t<u$ and $z_t'=z_{u-}z_u^{-1}z_t$
for $t\geq u$. Jumps may be successively removed at several time points and the resulting process is
independent of the order at which these operations are performed.

\begin{lemm} \label{leztuniquetriple}
Let $z_t$ be a rcll process in $G$. Then there is at most one
quadruple $(b,A,\eta,\nu)$ such that the $(b,A,\eta,\nu)$-martingale
property holds for $z_t$.
\end{lemm}

\noindent {\bf Proof} \ By Lemma~\ref{leztfixedjumpeta}, $\nu$ is determined by
the process $z_t$. Let $J=\{t>0$; $\nu_t\neq 0\}=\{u_1,u_2,u_3,\ldots\}$.
We will write the martingale in (\ref{ztmartprop}) as
\begin{equation}
M_tf = f(z_t) - \int_0^t Hf(z_\cdot) - \sum_{u\leq t}\int_G[f(z_{u-}x)-f(z_{u-})]\nu_u(dx),
\label{Mtfztmart}
\end{equation}
where $\int_0^tHf(z_\cdot)$ is the sum of integrals $\int_0^t(\cdots)db_j(s)$, $\int_0^t(\cdots)dA_{jk}(s)$
and $\int_0^t\int_G(\cdots)\eta(ds,dx)$ with bounded integrands. Because $b_i(t)$, $A_{jk}(t)$
and $\eta(t,\cdot)$ are continuous in $t$, one can show that $\int_0^tHf(z_\cdot)$ is
a bounded continuous function of $z_\cdot$ on $D(G)$ under the Skorohod metric.

Let $z_t^1$ be the process $z_t$
when its jump at time $u\in J$ is removed. Then for $t\geq v\geq u$,
\[M_tf - M_vf = f(z_t) - f(z_v) - \int_v^t Hf(z_\cdot) - \sum_{v<s\leq t}\int_G [f(z_{s-}x)-
f(z_{s-})]\nu_s(dx).\]
Replacing $f$ by $f'=f\circ(z_{u-}z_u^{-1})$ in the above yields
\[M_tf' - M_vf' = f(z_t^1) - f(z_v^1) - \int_v^t Hf(z_\cdot^1) - \sum_{v<s\leq t}
[f(z_{s-}^1x)-f(z_{s-}^1)]\nu_s(dx).\]
Because $E[M_tf'-M_vf'\mid\myF_v^z]=0$, it is then easy to show that $M_tf$
in (\ref{Mtfztmart}) is still a martingale under $\myF_t^z$ when $z_t$ and $J$ are replaced
by $z_t^1$ and $J-\{u\}$ respectively.

The jumps of $z_t$ at $u_1,u_2,\ldots,u_n$ may be successively removed
to obtain a new process $z_t^n$ such that $M_tf$ is still
a martingale with $z_t$ and $J$ replaced by $z_t^n$ and $J^n=J-\{u_1,u_2,\ldots,u_n\}$,
that is, $z_t^n$ has the $(b,A,\eta,\nu^n)$-martingale property
with $\nu^n=\{\nu_u$; $u\in J^n\}$.

We will show that along a subsequence of $n\to\infty$, $z_t^n$ converges weakly under the Skorohod metric to a rcll process $z_t'$ such
that $M_t'f$, given by (\ref{ztmartprop}) with $z_t$ and $J$ replaced by $z_t'$ and $\emptyset$, is a martingale.
This means that $z_t'$ has the $(b,A,\eta,0)$-martingale property, and then by Lemma~\ref{lefein}, $(b,A,\eta)$ is determined by $z_t'$,
which will prove Lemma~\ref{leztuniquetriple}.

We will start with a computation similar to the one in Section~\ref{secAB}.
Fix $T>0$ and a neighborhood $U$ of $e$. Let $\tau$ be the $U^c$-displacement time
of process $z_t^n$ from a stopping time $\sigma\leq T$. For $f\in C_c^\infty(G)$
with $f(e)=1$, $f=0$ on $U^c$ and $0\leq f\leq 1$ on $G$, let $f_\sigma(z)=f((z_\sigma^n)^{-1}z)$.
Then for $\delta>0$,
\begin{eqnarray}
&& P(\tau<\delta) = E[f_\sigma(z_\sigma^n)-
f_\sigma(z_{\sigma+\tau}^n);\,\tau<\delta] \leq
E[f_\sigma(z_\sigma^n) - f_\sigma(z_{\sigma+\tau\wedge\delta}^n)] \nonumber \\
&=& E[M_\sigma^nf_\sigma - M_{\sigma+\tau\wedge\delta}^nf_\sigma] -
E[\int_\sigma^{\sigma+\tau\wedge\delta}Hf_\sigma(z_\cdot^n)] \nonumber \\
&& - E\{\sum_{\sigma<u\leq\sigma+\tau\wedge\delta,\,u\in J^n}\int_G[f_\sigma(z_{u-}^nx)-f_\sigma(z_{u-}^n)]\nu_u(dx)\}
\nonumber \\
&\leq& 0 + \sum_j E\{\int_\sigma^{\sigma+\delta} \vert\cdots\vert db_j(s)\}\ + \sum_{j,k}
E\{\int_\sigma^{\sigma+\delta} \vert\cdots\vert dA_{jk}(s)\} \nonumber \\
&& \hspace{-0.3in} + \ E\{\int_\sigma^{\sigma+\delta}\int_G \vert\cdots\vert\eta(ds,dx)\} +
E\{\sum_{\sigma<u\leq\sigma+\delta,\,u\in J^n}\vert\int_G[f_\sigma(z_{u-}^nx)-f(z_{u-}^n)]\nu_u(dx)\vert\}.
\label{Ptaudelta2}
\end{eqnarray}
Because $b_j(t)$, $A_{jk}(t)$ and $\eta(t,\cdot)$ are continuous in $t$, and the associated integrands
are bounded, the first three nonzero terms in (\ref{Ptaudelta2}) converge to $0$ as $\delta\to 0$ uniformly in $n$
and $\sigma\leq T$. By Taylor's expansion, $f(z_{u-}^nx)-f(z_{u-}^n)=\sum_j\xi_jf(z_{u-}^n)\phi_j(x)+
O(\Vert\phi_\cdot\Vert^2)$ for $x$ in a neighborhood $U$ of $e$. By (\ref{sumnu2}),
it can be shown that the last term in (\ref{Ptaudelta2}) converges to $0$ as $n\to\infty$
uniformly in $\delta$ and $\sigma\leq T$. It now follows that $\lims_{n\to\infty}\sup_{\sigma\leq T}
P(\tau<\delta)\to 0$ as $\delta\to 0$, which verifies condition (\ref{supPtauU}) in
Lemma~\ref{lerelweakcompt}. Because the jumps of $z_t^n$ are those of $z_t$,
\[\sup_n\sup_{\sigma\leq T}P[(z_{\sigma-}^n)^{-1}z_\sigma^n\in K^c]
\ \leq\ \sup_{\sigma\leq T}P[z_{\sigma-}^{-1}z_\sigma\in K^c]\ \leq
\ P[z_{t-}^{-1}z_t\in K^c\ \mbox{for some $t\leq T$}]\ \downarrow 0\]
  as compact $K\uparrow G$ by the rcll property of $z_t$. This
verifies condition (\ref{supPxsigma}). Now
Lemma~\ref{lerelweakcompt} may be applied to show the weak convergence of $z_t^n$ in Skorohod metric to a rcll process $z_t'$ in $G$.

Because $\int_0^t Hf(z_\cdot)$ is a bounded continuous function on $D(G)$,
for any bounded continuous function $F(z_\cdot)$ on $D(G)$, $E\{[\int_0^t
Hf(z_\cdot^n)]F(z_\cdot^n)\}\to E\{\int_0^t Hf(z_\cdot')]
F(z_\cdot')\}$ as $n\to\infty$. In particular, this holds when $F$ is measureable
under $\sigma\{z_u$; $u\leq s\}$ for $s<t$. Because $M_tf$ in (\ref{Mtfztmart}) is a martingale when $z_t$ and $J$
are replaced by $z_t^n$ and $J^n$, and the sum $\sum_{u\leq t}$ in (\ref{Mtfztmart}) with $J$
replaced by $J^n$ converges to $0$ as $J^n\to\emptyset$, it follows that $M_tf$ is still a martingale
when $z_t$ and $J$ are replaced by $z_t'$ and $\emptyset$. \ $\Box$
\vspace{2ex}

To prove the uniqueness of the extended \levy triple $(b,A,\eta)$ in Theorem~\ref{thmartrep},
we will need a transformation rule for the martingale property (\ref{ztmartprop}) in the lemma below.

\begin{lemm} \label{lemartproptrans}
Let $z_t$ be a rcll process in $G$ having the $(b,A,\eta,\nu)$-martingale property.
If $u_t$ is a drift of finite variation with components $u_i(t)$, then $z_tu_t$ has
the $(b^u,A^u,\eta^u,\nu^u)$-martingale property, where $dA^u(t)=
[\Ad(u_t^{-1})]dA(t)[\Ad(u_t^{-1})]'$, $\eta^u(dt,f)=\int_G
f(u_t^{-1}xu_t)\eta(dt,dx)$ and $\nu_t^u(f)=\int_G
f(u_{t-}^{-1}xu_{t-})\nu_t(dx)$ for $f\in\myB_+(G)$, and $b_t^u$ is
given in components by
\[db_i^u(t) = \sum_p[\Ad(u_t^{-1})]_{ip}db_p(t) + du_i(t) +
\int_G\{\phi_i(u_t^{-1}xu_t) - \sum_p[\Ad(u_t^{-1})]_{ip}\phi_p(x)\}
\eta(dt,dx).\]
\end{lemm}

\noindent {\bf Proof} \ The proof is similar to the proof of
Lemma~\ref{lemartfinvar}, approximating $u_t$ by step functions $u_t^n=u_{t_{ni}}$, $t_{ni}\leq
t<t_{n\,i+1}$ and applying the $(b,A,\eta,\nu)$-martingale property for $z_t$ with the
function $f_{ni}(z)=f(zu_{t_{ni}})$, $f\in C_c^\infty(G)$.
We omit the essentially repetitive details. \ $\Box$

\begin{lemm} \label{leuniquetriple}
There is at most one extended \levy triple $(b,A,\eta)$ which
represents a given inhomogeneous \levy process $x_t$ in $G$ as in Theorem~\ref{thmartrep}.
\end{lemm}

\noindent {\bf Proof} \ Suppose two extended \levy triples $(b^1,A^1,\eta^1)$ and $(b^2,A^2,\eta^2)$
represent the same inhomogeneous \levy process $x_t$ in $G$. By Corollary~\ref{coxtfixedjumpeta}, $\eta^1$
and $\eta^2$ have the same discontinuous part, and hence $b_t^1$ and $b_t^2$ have the same jumps.
Let $x_t=z_t^1b_t^1=z_t^2b_t^2$. Then $z_t^1=z_t^2u_t$, where $u_t=b_t^2(b_t^1)^{-1}$ is a rcll path in $G$.
Note that $u_t=(z_t^2)^{-1}z_t^1$ is a semimartingale, that is, $f(u_t)$ is a real semimartingale
for any $f\in C_c^\infty(G)$. Because $u_t$ is non-random, it then follows that $u_t$ has a finite variation.
We now show $u_t$ is continuous and so is a drift. Let $x=u_{t-}^{-1}u_t=b_{t-}^1(b_{t-}^2)^{-1}b_t^2
(b_t^1)^{-1}$, then $(b_t^1)^{-1}xb_t^1=
[(b_{t-}^1)^{-1}b_t^1]^{-1}[(b_{t-}^2)^{-1}b_t^2]=e$ because $b_t^1$
and $b_t^2$ have the same jumps. This implies $x=e$ and hence $u_t$ is continuous.

By Lemma~\ref{leztmartproperty}, $z_t^i$ has the $(\bar{b}^i,\bar{A}^i,\bar{\eta}^i,
\bar{\nu}^i)$-martingale property for $i=1,2$. By Lemma~\ref{lemartproptrans}, $z_t^1$
also has the $(\bar{b}^{2\,u},\bar{A}^{2\,u},\bar{\eta}^{2\,u},\bar{\nu}^{2\,u})$-martingale
property. By Lemma~\ref{leztuniquetriple}, $\bar{b}^1=\bar{b}^{2\,u}$, $\bar{A}^1=
\bar{A}^{2\,u}$ and $\bar{\eta}^1=\bar{\eta}^{2\,u}$. Then for any $f\in C_c^\infty(G)$,
\begin{eqnarray*}
&& \int_G f(b_t^1x(b_t^1)^{-1})\eta^1{}^c(dt,dx) =
\bar{\eta}^1(dt,f) = \bar{\eta}^{2\,u}(dt,f) =
\int_G f(u_t^{-1}xu_t)\bar{\eta}^2(dt,dx)  \\
&=& \int_G f(u_t^{-1}b_t^2x(b_t^2)^{-1}u_t)\eta^2{}^c(dt,dx) = \int_G
f(b_t^1x(b_t^1)^{-1})\eta^2{}^c(dt,dx).
\end{eqnarray*}
This implies $\eta^1{}^c=\eta^2{}^c$ and hence $\eta^1=\eta^2$. Moreover,
\begin{eqnarray*}
  d\bar{b}_i^{2\,u}(t) &=& \sum_j[Ad(u_t^{-1})]_{ij}d\bar{b}_j^2(t) +
  du_i(t) + \int_G\{\phi_i(u_t^{-1}xu_t) - \sum_j\phi_j(x)
[\Ad(u_t^{-1})]_{ij}\}\bar{\eta}^2(dt,dx) \\
&=& \int_G\{\sum_j[\Ad(u_t^{-1})]_{ij}\phi_j(b_t^2x(b_t^2)^{-1}) -
\sum_p\phi_p(x)[\Ad(u_t^{-1}b_t^2)]_{ip}\}\eta^2{}^c(dt,dx) + du_i(t) \\
&&\ \ + \int_G\{\phi_i(u_t^{-1}b_t^2x(b_t^2)^{-1}u_t) -
\sum_j\phi_j(b_t^2x
(b_t^2)^{-1})[\Ad(u_t^{-1})]_{ij}\}\eta^2{}^c(dt,dx) \\
&=& \int_G\{\phi_i(u_t^{-1}b_t^2x(b_t^2)^{-1}u_t) - \sum_p\phi_p(x)
[\Ad(u_t^{-1}b_t^2)]_{ip}\}\eta^2{}^c(dt,dx) + du_i(t) \\
&=& \int_G\{\phi_i(b_t^1x(b_t^1)^{-1}) - \sum_p\phi_p(x)
[\Ad(b_t^1)]_{ip}\}\eta^1{}^c(dt,dx) + du_i(t) = d\bar{b}_i^1(t)\ +\
du_i(t).
\end{eqnarray*}
Because $\bar{b}^1=\bar{b}^{2\,u}$, this implies that $u_i(t)=0$ and hence $b^1=b^2$.
Now $\bar{A}^1=\bar{A}^{2\,u}=\bar{A}^2$, and with $b^1=b^2$, it follows that $A^1=A^2$.
The uniqueness of $(b,A,\eta)$ is proved. \ $\Box$

\section{Uniqueness and existence of the process} \label{secuniqueexistprocess}

\begin{lemm} \label{lerelweakcompt2}
Let $x_t^n$ be a sequence of rcll processes in $G$ as in Lemma~\ref{lerelweakcompt} and $u_n\in\bR_+$.
For each $n$, let $x_t^{n,m}$ be the process obtained from $x_t^n$ when possible fixed jumps at times $u_1,u_2,\ldots,u_m>0$
are removed. If for any $T>0$, $\eta>0$ and neighborhood $U$ of $e$, there are $\delta>0$ and an integer $m\geq 1$ such that
\begin{equation}
\overline{\lim_{n\to\infty}}\sup_{\sigma\leq
T}P(\tau_U^{\sigma,n,m}<\delta) \leq \eta, \label{supPtauUm}
\end{equation}
where $\sup_{\sigma\leq T}$ is taken over all stopping times
$\sigma\leq T$ and $\tau_U^{\sigma,n,m}$ is $\tau_U^\sigma$ for
process $x_t^{n,m}$, and if (\ref{supPxsigma}) also holds, then
a subsequence of $x_t^n$ converges weakly in $D(G)$.
\end{lemm}

\noindent {\bf Proof} \ Using the notation in the proof of Lemma~\ref{lerelweakcompt},
let $\tau_i^{\varepsilon,n,m}$ be the $\varepsilon$-displacement time $\tau_i^\varepsilon$ for process $x_t^{n,m}$. By (\ref{Fdelta})
and (\ref{supPtauUm}), for any $\eta$ and $\varepsilon>0$, there are $\delta$ and $m$ such that
\[\lims_{n\to\infty}P[\min_{i\geq 0}\{\tau_{i+1}^{\varepsilon,n,m}-\tau_i^{\varepsilon,n,m};
\,\tau_i^{\varepsilon,n,m}<T\}\leq\delta]\ \leq\ \eta.\]
Then by (\ref{mintauomega}), $\lims_{n\to\infty} P[\omega'(x^{n,m},\delta,T)>2\varepsilon]\leq\eta$.

The computation of $\omega'(x,\delta,T)$ in (\ref{omegaprime}) is based on the oscillations over partitions that cover $[0,\,T]$
with spacing $>\delta$, called $\delta$-partitions. Let $J_m=\{u_1,\ldots,u_m\}$. We may assume $\delta\ <\ \min\{\vert u\vert\wedge\vert u-v\vert;
\ u,v\in J_m\ {\rm with}\ u\neq v\}$, and either $T\in J_m$ or $\vert T-u\vert>\delta$ for $u\in J_m$. Suppose there is a $\delta$-partition
with oscillation $\leq 2\varepsilon$. Then any interval of this $\delta$-partition contains at most one point in $J_m$.
Because adding more partition points will not increase oscillation,
adding $J_m$ to the $\delta$-partition together with the midpoints of intervals
that do not intercept $J_m$, and then suitably combining intervals,
we obtain a $(\delta/2)$-partition containing $J_m$, with oscillation $\leq 4\varepsilon$.
Then $\lims_{n\to\infty}P[\omega'(x^n,\delta/2,T)> 4\varepsilon]\leq\eta$. This verifies (\ref{omegapvarepsilon}).

It remains to verify (\ref{xtnKc}). If $\omega'(x,\delta,T)<\varepsilon$,
then there is a $\delta$-partition $\{t_i\}$ with oscillation $<\varepsilon$.
This implies that $\tau_1^\varepsilon\not\in[0,\,t_1)$ and $\tau_i^\varepsilon$ belong
to different intervals of the partition. Thus
\begin{equation}
\omega'(x,\delta,T)<\varepsilon\ \ \Longrightarrow\ \ \min_{i\geq
0}\{\tau_{i+1}^\varepsilon-\tau_{i-1}^\varepsilon;\,\tau_i^\varepsilon<T\}
>\delta\ \ \ \ \mbox{(setting $\tau_{-1}^\varepsilon=0$).}
\label{omegamintau}
\end{equation}
By (\ref{omegapvarepsilon}) and (\ref{omegamintau}), $\lims_{n\to\infty}P(A_n)\to 0$ as $\delta\to 0$,
where $A_n=[\min_{i\geq 0}\{\tau_{i+1}^{\varepsilon,n}-
\tau_{i-1}^{\varepsilon,n}; \,\tau_i^{\varepsilon,n}<T\}<\delta]$.
The rest of proof is very similar to the last part of the proof of
Lemma~\ref{lerelweakcompt}. \ $\Box$
\vspace{2ex}

Let $z_t$ be an inhomogeneous \levy process in $G$ having the $(b,A,\eta,\nu)$-martingale property
and let $J=\{u>0$; $\nu_u\neq 0\}=\{u_1,u_2,\ldots\}$.
For $p<q$, let $J_p^q=\{u_1,\ldots,u_p,u_{q+1}$, $u_{q+2},\ldots\}$ and let $z_t^{p,q}$
be the process $z_t$ when its fixed jumps at $u_{p+1},\ldots,u_q$ are removed. As in the proof
of Lemma~\ref{leztuniquetriple}, it can be shown that $z_t^{p,q}$ has the $(b,A,\eta,\nu^{p,q})$-martingale
property with $\nu^{p,q}=\{\nu_u$; $u\in J_p^q\}$. Let $J_p=\{u_1,u_2,\ldots,u_p\}$.

\begin{lemm} \label{leztpq}
The family of processes, $z_t^{p,q}$ for $p<q$, are weakly compact in $D(G)$. Moreover,
along a subsequence of $n\to\infty$, $z_t^{p,n}$ converges weakly to an inhomogeneous \levy
process $z_t'^p$ that has the $(b,A,\eta,\nu'^p)$-martingale property with $\nu'^p=\{\nu_u$; $u\in J_p\}$.
\end{lemm}

\noindent {\bf Proof} \ Let $z_t^n=z_t^{p_n,q_n}$, where $p_n<q_n$
and $q_n\to\infty$ as $n\to\infty$. As in the proof of Lemma~\ref{leztuniquetriple},
but using Lemma~\ref{lerelweakcompt2} instead of Lemma~\ref{lerelweakcompt}, it can be shown
that a subsequence of $z_t^n$ converges weakly in $D(G)$ to some $z_t'$. Note that
the last term in (\ref{Ptaudelta2}) now takes the form $\sum_{\sigma<u\leq\sigma+\delta}(\cdots)$
with $u\in\{u_1,\ldots,u_{p_n},u_{q_n+1},\ldots\}$, which may not tend to $0$ as $n\to\infty$,
but we may remove finitely many fixed jumps to make it arbitrarily small. By Lemma~\ref{lerelweakcompt2},
one can show that the family of the processes $z_t^{p,q}$, $p<q$, are weakly compact on $D(G)$.

Let $z_t'^p$ be the weak limit of $z_t^{p,n}$ as $n\to\infty$. It is easy to see
that $z_t'^p$ is an inhomogeneous \levy process.
We now prove that $z_t'^p$ has the $(b,A,\eta,\nu'^p)$-martingale property,
that is, $M_tf$ given in (\ref{Mtfztmart}) is a martingale
when $z_t$ and $J$ are replaced by $z_t'^p$ and $J_p$. As mentioned in the proof
of Lemma~\ref{leztuniquetriple}, $\int_0^t Hf(z_\cdot)$ in (\ref{Mtfztmart}) is a bounded continuous function
on $D(G)$. It suffices to show that for any $u\in J_p$ and any bounded continuous function $F(z)$
on $D(G)$, $\myF_v^z$-measurable for some $v<u$, as $n\to\infty$,
\begin{equation}
E\{F(z_\cdot^{p,n})\int[f(z_{u-}^{p,n}x)-f(z_{u-}^{p,n})]\nu_u(dx)\}\ \to\ E\{F(z_\cdot'^p)
\int[f(z_{u-}'^px)-f(z_{u-}'^p)]\nu_u(dx)\}. \label{EFintnu}
\end{equation}

For $\varepsilon>0$, choose $\phi_\varepsilon\in C(\bR)$ such that $0\leq\phi_\varepsilon
\leq 1$, $\phi_\varepsilon(u)=1$ and $\phi_\varepsilon(t)=0$ for $\vert t-u\vert>\varepsilon$.
For any $g\in C(G)$ such that $0\leq g\leq 1$ and $g=0$ near $e$, let $H_\varepsilon(z_\cdot)=\sum_{t>0}\phi_\varepsilon(t)
g(z_{t-}^{-1}z_t)$. Then $\tilde{H}_\varepsilon(z_\cdot)=H_\varepsilon(z_\cdot)\wedge 2$ is a bounded continuous function on $D(G)$
and hence $E[\tilde{H}_\varepsilon(z_\cdot^{p,n})]\to E[\tilde{H}_\varepsilon
(z_\cdot'^p)]$ as $n\to\infty$. Because $\tilde{H}_\varepsilon(z)\to g(z_{u-}^{-1}z_u)$ as $\varepsilon\to 0$, $\tilde{H}_\varepsilon(z)
\leq\sum_{t\leq T,\vert t-u\vert\leq\varepsilon}g(z_{t-}^{-1}z_t)$ and $E[\sum_{t\leq T}g(z_{t-}^{-1}z_t)]\leq\sum_{t\leq T}\nu_t(g)<\infty$,
it follows that $E[\tilde{H}_\varepsilon(z_\cdot^{p,n})]\to\nu_u(g)$ as $\varepsilon\to 0$ uniformly in $n$. This implies
that $E[g((z_{u-}'^p)^{-1}z_u'^p)]=\lim_{\varepsilon\to 0}E[\tilde{H}_\varepsilon(z_\cdot'^p)]=\nu_u(g)$ and hence the distribution
of $(z_{u-}'^p)^{-1}z_u'^p$ is $\nu_u$.

We may assume $F$ and $f$ in (\ref{EFintnu}) are bounded by $1$ in absolute values. Let $J_\varepsilon(z_\cdot)=F(z_\cdot)\sum_{t>0}\phi_\varepsilon(t)
[f(z_t)-f(z_{t-})]g(z_{t-}^{-1}z_t)$ and let $\tilde{J}_\varepsilon(z_\cdot)=[J_\varepsilon(z_\cdot)\wedge 3]\vee(-3)$, where $a\vee b=\max(a,b)$.
Then $\tilde{J}_\varepsilon(z_\cdot)$ is a bounded continuous function on $D(G)$ and hence $E[\tilde{J}_\varepsilon(z_\cdot^{p,n})]\to
E[\tilde{J}_\varepsilon(z_\cdot'^p)]$ as $n\to\infty$. Because $\tilde{J}_\varepsilon(z_\cdot)\to F(z_\cdot)[f(z_u)-f(z_{u-})]g(z_{u-}^{-1}z_u)$ as $\varepsilon
\to 0$ and $\vert\tilde{J}_\varepsilon(z_\cdot)\vert\leq 2\sum_{t\leq T,\vert t-u\vert\leq\varepsilon}g(z_{t-}^{-1}z_t)$, it follows
that $E[\tilde{J}_\varepsilon(z_\cdot^{p,n})]\to E\{F(z_\cdot^{p,n})[f(z_u^{p,n})-f(z_{u-}^{p,n})]g((z_{u-}^{p,n})^{-1}z_u^{p,n})\}$ as $\varepsilon
\to 0$ uniformly in $n$. Then $E\{F(z_\cdot'^p)[f(z_u'^p)-f(z_{u-}'^p)]g((z_{u-}'^p)^{-1}z_u'^p)\}=\lim_{\varepsilon\to 0}E[\tilde{J}_\varepsilon(z_\cdot'^p)]=
\lim_{\varepsilon\to 0}\lim_{n\to\infty}E[\tilde{J}_\varepsilon(z_\cdot^{p,n})]=\lim_{n\to\infty}\lim_{\varepsilon\to 0}
E[\tilde{J}_\varepsilon(z_\cdot^{p,n})]$. Letting $g\uparrow 1_{G-\{e\}}$ proves (\ref{EFintnu}). \ $\Box$

\begin{lemm} \label{leuniqueindist}
The distribution of an inhomogeneous \levy process $z_t$ in $G$ with $z_0=e$ having
the martingale property under a given quadruple $(b,A,\eta,\nu)$ is unique.
\end{lemm}

\noindent {\bf Proof} \ Because $z_t'^p$ in Lemma~\ref{leztpq} has only finitely many fixed jumps, after
removing these fixed jumps, it becomes stochastically continuous. By Lemma~\ref{lefein}, the distribution
of $z_t'^p$ is completely determined by $(b,A,\eta,\nu'^p)$, and hence also by $(b,A,\eta,\nu)$.

By Lemma~\ref{leztpq}, the processes $x_t^{p,q}$, $p<q$, are weakly compact on $D(G)$.
Then by Remark~7.3 in \cite[chapter~3]{ek}, for any $\varepsilon>0$,
there is a compact $K_\varepsilon\subset G$ such that
\begin{equation}
\forall \ p<q,\ \ \ \ P(z_s^{p,q}\in K_\varepsilon\ {\rm for}\ s\leq t)
\ \geq\ 1-\varepsilon. \label{PzspqK}
\end{equation}
Let $f\in C_c^\infty(G)$ with $\vert f\vert\leq 1$. We will show
that $\forall\varepsilon>0$, $\exists$ constant $c_\varepsilon=c_{\varepsilon,f,t}>0$ such that
\begin{equation}
\forall \ p<q,\ \ \ \ \vert E[f(z_t)] - E[f(z_t^{p,q})]\vert\ \leq\ c_\varepsilon\sum_{i=p+1}^q[\sum_{j=1}^d
\vert\nu_{u_i}(\phi_j)\vert + \nu_{u_i}(\Vert\phi_\cdot\Vert^2)+\nu_{u_i}(U^c)] + 2\varepsilon \label{Efztfztpq}
\end{equation}
for a small neighborhood $U$ of $e$.

Recall $z_t^{p,q}$ is obtained from $z_t$ after the fixed jumps at times $u_{p+1},u_{p+2},\ldots,u_q$ are
removed. Let $v_1<v_2<\cdots<v_{q-p}$ be the ordered values of these time points and let $m\leq q-p$ be the largest
integer such that $v_m\leq t$. Let $\sigma_i=z_{v_i-}^{-1}z_{v_i}$ and for $u<v$, let $z_{u,v}=z_u^{-1}z_{v-}$
and $z_{u,v}^{p,q}=(z_u^{p,q})^{-1}z_{v-}^{p,q}$.
For $x_i\in G$, let $\prod_{i=1}^nx_i=x_1x_2\cdots x_n$. Set $v_0=0$ and $z_{u,u}=
z_{u,u}^{p,q}=e$, and let $H=\cap_{i=1}^m[z_{0,v_i}^{p,q}\in K_\varepsilon]$. By (\ref{PzspqK}), $P(H)\geq 1-\varepsilon$.
For simplicity, assume $z_t=z_{t-}$.
\begin{eqnarray*}
&& f(z_t)-f(z_t^{p,q}) = f(z_{0,v_1}^{p,q}\sigma_1z_{v_1,t}) - f(z_{0,v_m}^{p,q}z_{v_m,t}) = \sum_{i=1}^m
[f(z_{0,v_i}^{p,q}\sigma_iz_{v_i,t}) - f(z_{0,v_i}^{p,q}z_{v_i,t})] \\
&=& \sum_{i=1}^m[f(z'\sigma z) - f(z'z)]\ \ \ \ \mbox{($z'=z_{0,v_i}^{p,q}$, $z=z_{v_i,t}$ and $\sigma=\sigma_i$)} \\
&=& \sum_{i=1}^m[f(z'e^{\sum_j\phi_j(\sigma)\xi_j}z)-f(z'z)]1_{[\sigma\in U]} + \sum_{i=1}^m[f(z'\sigma z)-f(z'z)]1_{[\sigma\in U^c]} \\
&=& \sum_{i=1}^m[f(e^{\sum_{j,k}[\Ad(z')]_{kj}\phi_j(\sigma)\xi_k}z'z) - f(z'z)]1_{[\sigma\in U]} + \sum_{i=1}^m
[f(z'\sigma z)-f(z'z)]1_{[\sigma\in U^c]} \\
&=& \sum_{i=1}^m\{\sum_{j,k=1}^d \xi_k^rf(z'z) [\Ad(z')]_{kj}\phi_j(\sigma) + O(\Vert\phi_\cdot(\sigma)
\Vert^2)\}1_{[\sigma\in U]} + \sum_{i=1}^m [f(z'\sigma z)-f(z'z)]1_{[\sigma\in U^c]}
\end{eqnarray*}
by Taylor's formula, where $\xi^rf(z)=\frac{d}{dt}
f(e^{t\xi}z)\mid_{t=0}$. Note that $z_t^{p,q}$ and $H$ are independent of $\sigma_i$, and on $H$, $z_s^{p,q}
\in K_\varepsilon$ and hence $[\Ad(z_s^{p,q})]_{kj}$ are bounded for all $s\leq t$. It follows that
\begin{eqnarray*}
&& \vert E[f(z_t)-f(z_t^{p,q})]\vert\ \leq\ \vert E\{[f(z_t)-f(z_t^{p,q})]1_H\}\vert + 2\varepsilon \\
&\leq& c_\varepsilon'\sum_{i=1}^m[\sum_j\vert\nu_{v_i}(\phi_j1_U)\vert + O(\Vert\phi_\cdot\Vert^2) + \nu_{v_i}(U^c)] + 2\varepsilon
\end{eqnarray*}
for some constant $c_\varepsilon'>0$. This proves (\ref{Efztfztpq}) because $\vert\nu_u(\phi_j1_U)\vert
\leq\vert\nu_u(\phi_j)\vert+(\sup_G\vert\phi_j\vert)\nu_u(U^c)$.

By (\ref{sumnu2}), letting first $q\to\infty$ and then $p\to\infty$ in (\ref{Efztfztpq}) yields $E[f(z_t)]-E[f(z_t'^p)]\to 0$.
Because $E[f(z_t'^p)]$ is determined by $(b,A,\eta,\nu)$, so is $E[f(z_t)]$. This proves the uniqueness of
the one-dimensional distribution of an inhomogeneous \levy process $z_t$ having the $(b,A,\eta,\nu)$-martingale property.
To complete the proof, one just need to use the standard argument to derive the uniqueness of distribution from
the uniqueness of one-dimensional distribution and the martingale property, see the proof
of Theorem~4.2(a) in \cite[chapter~4]{ek}. Note the setups here and in \cite{ek} are not quite the same,
but the argument is essentially the same. \ $\Box$
\vspace{2ex}

By Lemmas \ref{leztmartproperty} and \ref{leuniqueindist},
we obtain the unique distribution for the process $x_t$ in Theorem~\ref{thmartrep}.

\begin{coro} \label{couniquedist}
The distribution of an inhomogeneous \levy process $x_t$ in $G$ with $x_0=e$
represented by a given extended \levy triple is unique.
\end{coro}

\noindent {\bf Proof of Theorem~\ref{thmartrep}} \ It remains to
prove that given an extended \levy triple $(b,A,\eta)$, there is an
inhomogeneous \levy process $x_t$ in $G$ represented by $(b,A,\eta)$,
that is, with $x_t=z_tb_t$, (\ref{Mtf}) is a martingale under
the natural filtration $\myF_t^x$ of $x_t$ for any $f\in C_c^\infty(G)$.

Let $(\bar{b},\bar{A},\bar{\eta},\bar{\nu})$ be as in
Lemma~\ref{leztmartproperty}. Then $(\bar{b},\bar{A},\bar{\eta})$ is
a \levy triple with $\bar{b}_t$ of finite variation. By
Lemma~\ref{lefein}, there is a stochastically continuous inhomogeneous
\levy process $z_t'$ having the $(\bar{b},\bar{A},\bar{\eta},0)$-martingale property.
This means that for $f\in C_c^\infty(G)$, $f(z_t')-\int_0^t\bar{H}f(z_\cdot')$ is a martingale
under the natural filtration of $z_t'$, where $\int_0^t\bar{H}f(z_\cdot')$ is
the expression $\int_0^t Hf(z_\cdot)$ in (\ref{Mtfztmart}) with $b,A,\eta,z_\cdot$ replaced
by $\bar{b},\bar{A},\bar{\eta},z_\cdot'$.

Let $J=\{u_1,u_2,u_3,\ldots\}$ be the set of discontinuity time
points of $\eta$, and let $\{v_u$, $u\in J\}$ be a family of $G$-valued random variables
with distributions $\nu_u$ such that they are mutually independent
and independent of process $z_t'$. For $u=u_1$, and let $z_t^1$ be the
process $z_t'$ when the fixed jump $\bar{v}_u=b_{u-}v_uh_u^{-1}
b_{u-}^{-1}$ is added at time $u$, that is, $z_t^1=z_t'$ for $t<u$
and $z_t^1=z_{u-}'\bar{v}_uz_u'^{-1}z_t'$ for $t\geq u$. Then it is
easy to show that $f(z_t^1)-\int_0^t\bar{H}f(z_\cdot^1)-1_{[u\leq
t]}\int_G[f(z_{u-}^1x)- f(z_{u-}^1)]\bar{\nu}_u(dx)$ is a martingale under
the natural filtration of $z_t^1$. More generally, let $z_t^n$ be the process $z_t'$
when the fixed jumps $\bar{v}_{u_1},\bar{v}_{u_2},\ldots,\bar{v}_{u_n}$ are successively
added at times $u_1,u_2,\ldots,u_n$. Then
\begin{equation}
f(z_t^n) - \int_0^t \bar{H}f(z_s^n)d\bar{b}(s) - \sum_{u\leq t,\,u\in J_n}\int_G[f(z_{u-}^nx)-f(z_{u-}^n)]
\bar{\nu}_u(dx) \label{fztprimenu}
\end{equation}
  is a martingale under the natural filtration of $z_t^n$, where $J_n=\{u_1,u_2,\ldots,u_n\}$.

  As in the proof of Lemma~\ref{leztpq}, we may use Lemma~\ref{lerelweakcompt2} to show that a subsequence
of $z_t^n$ converge weakly in $D(G)$ to a rcll process $z_t$, and (\ref{fztprimenu}) is still
a martingale when $z_t^n$ and $J_n$ are replaced by $z_t$ and $J$.
Therefore, with $x_t=z_tb_t$, (\ref{Mtf}) is a martingale under $\myF_t^z=\myF_t^x$.

It remains to show that $x_t$ is an inhomogeneous \levy process in $G$.
Because $z_t^n$ is obtained by adding $n$ independent fixed
jumps to the inhomogeneous \levy process $z_t'$, so $z_t^n$ is an
inhomogeneous \levy process. The same is true for $z_t$, as the
weak limit of $z_t^n$, and hence true also for $x_t=z_tb_t$.  \ $\Box$

\section{Results on homogeneous spaces} \label{sec6}

Let $K$ be a closed subgroup of a Lie group $G$. The space $G/K$
of left cosets $gK$, $g\in G$, is called a homogeneous space. It is
equipped with the unique manifold structure under which the natural $G$-action
on $G/K$, $(g,g'K)\mapsto gg'K$ for $g,g'\in G$, is smooth (Theorem~4.2 in \cite[chapter~II]{helgason1}.
Moreover, the natural projection $\pi$: $G\to G/K$, $g\mapsto gK$, is smooth and open.

Let $x_t$ be a Markov process $x_t$ in $G/K$ with rcll paths. It is called $G$-invariant if its transition function $P_t$ is $G$-invariant
in the sense that
\[\forall t\geq 0,\  g\in G\ {\rm and}\ f\in\myB_+(X),\ \ \ \ P_t(f\circ g) = (P_tf)\circ g.\]
When $K=\{e\}$, such a process becomes a left invariant Markov process in $G$, that is, a \levy process in $G$. Therefore, a rcll $G$-invariant Markov process
in $G/K$ will also be called a \levy process. More generally, a rcll inhomogeneous Markov process in $G/K$ with a $G$-invariant transition function $P_{s,t}$,
will be called an inhomogeneous \levy process.

Let $X$ be a manifold under the transitive action of a Lie group $G$ and let $K=\{g\in G$; $go=o\}$ be the fix-point subgroup at a point $o\in X$.
Then $X$ may be identified with the homogeneous space $G/K$ via the map $gK\mapsto go$ (Theorem~3.2 and Proposition~4.3
in \cite[chapter~II]{helgason1}), under which the $G$-action on $X$ is identified with the natural $G$-action on $G/K$. Although this identification is not
unique, as a different $o$ may be chosen, but (inhomogeneous) \levy processes in $X$, defined as $G$-invariant (inhomogeneous) Markov processes,
are independent of this identification.

In the rest of this paper, let $x_t$ be an inhomogeneous \levy process in $X=G/K$ and let $J$ be the countable set of fixed jump times of $x_t$.
We will assume $K$ is compact.

A measure $\mu$ on $X=G/K$ is called $K$-invariant if $k\mu=\mu$ for $k\in K$, where $k\mu$ is the measure defined
by $k\mu(f)=\mu(f\circ k)$ for $f\in\myB_+(X)$. Because $P_{s,t}$ is $G$-invariant, $\mu_{s,t} = P_{s,t}(o,\cdot)$ is a $K$-invariant
probability measure on $X$, where $o=eK$ is the origin in $X$.

A Borel measurable map $S$: $X\to G$ is called a section map if $\pi\circ S=\id_X$ (the identity map on $X$). Such a map is not unique and
it may not be continuous on $X$. However, one may always choose a section map that is smooth near any given point in $X$. Integrals like
\[\int f(xy)\mu(dy) = \int f(S(x)y)\mu(dy),\ \ \ \ \int f(xyz)\mu(dy)\nu(dz) =
\int f(S(S(x)y)z)\mu(dy)\nu(dz),\]
are well defined for $K$-invariant measures $\mu$ and $\nu$ on $X$, independent of the choice for the section map $S$. This is easily verified
by the $K$-invariance of $\mu$ and $\nu$, because if $S'$ is another section map, then for any $x\in X$, $S'(x)=S(x)k_x$ for some $k_x\in K$. Note also
that $S(S(x)y)=S(x)S(y)k$ for some $k=k(x,y)\in K$, so the second integral above may also be computed as $\int f(S(x)S(y)z)\mu(dy)\nu(dz)$.

A random variable in $X$ is called $K$-invariant if its distribution is so. It is clear that if $x$ and $y$ are two independent random variables in $X$
with $y$ being $K$-invariant, then the distribution of $xy=S(x)y$ is independent of the choice for section map $S$.

By the Markov property and the $G$-invariance of $P_{s,t}$,
\begin{eqnarray}
E[f(x_{t_1},x_{t_2},\ldots,x_{t_n})]\ &=&\ \int\mu_0(dx_0)
\mu_{0,t_1}(dx_1) \mu_{t_1,t_2}(dx_2)\cdots\mu_{t_{n-1},t_n}(dx_n)
\nonumber \\
&&\ \ f(x_0x_1,\,x_0x_1x_2,\,\ldots,\,x_0x_1x_2\cdots x_n)
\label{findimdist}
\end{eqnarray}
for $t_1<t_2<\cdots<t_n$ and $f\in\myB_+(G^n)$, where $\mu_0$ is the initial distribution.

\begin{prop} \label{prii}
A rcll process $x_t$ in $X$ is an inhomogeneous \levy process if and only
if it has independent increments in the sense that
for $s<t$, $x_s^{-1}x_t=S(x_s)^{-1}x_t$ is independent
of $\myF_s^x$ and has a $K$-invariant distribution $\mu_{s,t}$ independent of the choice
for the section map $S$.
\end{prop}

\noindent {\bf Proof} \ From (\ref{findimdist}),
it is easy to show that $S(x_s)^{-1}x_t$ is independent of $\myF_s^x$
and its distribution $\mu_{s,t}$ is independent of $S$. Conversely, assume
this property, then for $s<t$ and $f\in\myB_+(X)$, $E[f(x_t)\mid
\myF_s^x]=E[f(S(x_s)S(x_s)^{-1}x_t)\mid\myF_s^x]=\int_X f(S(x_s)y)\mu_{s,t}(dy)$.
If $\mu_{s,t}$ is $K$-invariant,
then $x_t$ is an inhomogeneous Markov process in $X$ with a $G$-invariant
transition function $P_{s,t}f(x)=\int_Xf(S(x)y)\mu_{s,t}(dy)$, and hence $x_t$ is
an inhomogeneous \levy process in $X$. \ $\Box$
\vspace{2ex}

The proof of Proposition~\ref{prii} may be slightly modified to show that a rcll process $x_t$ in $X$ is a \levy process if and only if it has independent
and stationary increments. Here the stationary increments mean the distribution of $x_s^{-1}x_t=S(x_s)^{-1}x_t$ depends only on $t-s$.

A measure function $\eta(t,\cdot)$, $t\geq 0$, on $X$ is defined
just as a measure function on $G$, that is, as a nondecreasing and
right continuous family of $\sigma$-finite measures on $X$.
It is called $K$-invariant if $\eta(t,\cdot)$ is $K$-invariant for each $t\geq 0$. As on $G$,
a measure function $\eta(t,\cdot)$ on $X$ may
be regarded as a $\sigma$-finite measure on $\bR_+\times X$ and it
may be written as a sum of its continuous part $\eta^c(t,\cdot)$
and its discontinuous part $\sum_{s\leq t}\eta(\{s\}\times\cdot)$ as in (\ref{etaetaprimeeta}).

Let $\Delta_n$: $0=t_{n0}<t_{n1}<t_{n2}<\cdots<t_{ni}
\uparrow\infty$ as $i\uparrow\infty$ be a sequence of partitions of $\bR_+$ with mesh $\Vert\Delta_n
\Vert\to 0$ as $n\to\infty$, and assume $\Delta_n$ contains $J$ as $n\to\infty$ as before. Let $x_t$ be an inhomogeneous L'{e}vy process in $X$ and
let $\mu_{ni}$ be the distribution of $x_{t_{n\,i-1}}^{-1}x_{t_{ni}}$ for $i\geq 1$. Define
\begin{equation}
\eta_n(t,\cdot) = \sum_{t_{ni}\leq t}\mu_{ni}, \label{etantX}
\end{equation}
  setting $\eta_n(0,\cdot)=0$. Then $\eta_n$ is a $K$-invariant measure function on $X$.

The proof of Proposition~\ref{pretanbdd} may be repeated on $X$,
regarding $x_s^{-1}x_t$ as $S(x_s)^{-1}x_t$ with the choice of a section map $S$.
Note that with a $G$-invariant metric $r$ on $X$, $r(x,y)=r(o,S^{-1}(x)y)$, and hence $S(x_s)^{-1}x_t\in U^c$
for a neighborhood $U$ of $o$ if and only if $r(x_s,x_t)>\delta$ for some $\delta>0$.
Then for any $t>0$ and any neighborhood $U$
of $o$,
\begin{equation}
\eta_n(t,U^c)\ \ {\rm is\ bounded\ in}\ n\ \ {\rm and}\ \
\eta_n(t,U^c)\downarrow 0\ \ {\rm as}\ U\uparrow X\ \ {\rm
uniformly\ in}\ n. \label{etanbddX}
\end{equation}

For $f\in C_b(X)$, let $\hat{f}=\int_K dk(f\circ k)$, where $dk$ is the normalized Haar measure
on $K$. Then $\hat{f}$ is $K$-invariant, that is, $\hat{f}\circ k=\hat{f}$
for $k\in K$, and for a $K$-invariant measure $\mu$ on $X$, $\mu(f)=\mu(\hat{f})$.
We note that to show the weak convergence of a sequence of $K$-invariant measures $\mu_n$ on $K$, it suffices
to show the convergence of $\mu_n(f)$ for $K$-invariant $f\in C_b(X)$.

Let $f\in C_b(X)$ be $K$-invariant and vanish near $o$. Then the function $(x,y)\mapsto f(S(x)^{-1}y)$ is independent
of the choice for the section map $S$ and is continuous. The proof of Proposition~\ref{pretantoeta} may be repeated
to show that as $n\to\infty$,
\[\eta_n(t,f) = \sum_{t_{ni}\leq t} E[f(S(x_{t_{n\,i-1}})^{-1}x_{t_{ni}})]
\ \to\ E[\sum_{s\leq t}f(S(x_{s-})^{-1}x_s)].\]

This shows that $\eta(t,\cdot)$ defined by $\eta(0,\cdot)=0$, $\eta(t,\{o\})=0$ and
\begin{equation}
\eta(t,f) = \eta(t,\hat{f}) = E[\sum_{s\leq t}\hat{f}(x_{s-}^{-1}x_s)]
\ \ \ \ \mbox{for $f\in C_b(G)$ vanishing near $o$}, \label{etatBX}
\end{equation}
is a well defined $K$-invariant measure function on $X$, independent of the choice of
the section map $S$ to represented $x_{s-}^{-1}x_s=S(x_{s-})^{-1}x_s$ in (\ref{etatBX}),
which will be called the jump intensity measure of process $x_t$. It is clear that $\eta=0$
if and only if the process $x_t$ is continuous, and $\eta(t,\cdot)$ is continuous in $t$
if and only if $x_t$ is stochastically continuous. Moreover, for any $t>0$ and $f\in C_b(X)$ vanishing
near $o$, $\eta_n(t,f)\to\eta(t,f)$ as $n\to\infty$.

A point $b\in X$ is called $K$-invariant if $kb=b$ for all $k\in K$.
For $x\in X$ and $K$-invariant $b$, the product $xb=S(x)b\in X$ is
well defined, independent of choice for the section map $S$.

\begin{prop} \label{prKinvpt}
Let $N$ be the set of $g\in G$ such that $go=\pi(g)$ are $K$-invariant.
Then $N$ is the normalizer of $K$, that is, $N=\{g\in G$; $gKg^{-1} = K\}$.
\end{prop}

\noindent {\bf Proof} \ Note that $g\in G$ with $go$ being $K$-invariant is characterized
by $g^{-1}Kg\subset K$. This implies that $g^{-1}Kg$ and $K$ have the same identity component $K_0=
g^{-1}K_0g$. Because $K$ is compact, its coset decomposition $K=K_0\cup[\cup_{i=1}^p(k_iK_0)]$ is
a finite union. Since $g^{-1}Kg=K_0\cup[\cup_{i=1}^p(g^{-1}k_igK_0)]$ is the coset decomposition
of $g^{-1}Kg$, it follows that $g^{-1}Kg=K$. \ $\Box$
\vspace{2ex}

By Proposition~\ref{prKinvpt}, the set of $K$-invariant points in $X$ has
a natural group structure with product $bb'=S(b)b'$ and
inverse $b^{-1}=S(b)^{-1}o$,  and the integral
\[\int f(xbyb')\mu(dy) = \int f(S(x)S(b)S(y)b')\mu(dy)\]
makes sense for $K$-invariant measure $\mu$ on $X$ and $K$-invariant
points $b,b'\in X$, independent of the choice of the section map $S$.

Let $\myg$ and $\myk$ be respectively the Lie algebras of $G$ and $K$.
Because $K$ is compact, there is a subspace $\myp$ of $\myg$ that is complementary
to $\myk$ and is $\Ad(K)$-invariant in the sense that $\Ad(k)\myp=\myp$
for $k\in K$. Choose a basis $\xi_1, \ldots,\xi_n$ of $\myp$.
There are $\phi_1,\phi_2,\ldots,\phi_n\in C_c^\infty(X)$ such
that $x=\exp[\sum_{i=1}^n\phi_i(x)\xi_i]o$ for $x$ near $o$. It then follows that
\begin{equation}
\forall k\in K,\ \ \ \ \sum_{i=1}^n \phi_i(x)[\Ad(k)\xi_i] = \sum_{i=1}^n\phi_i(kx)\xi_i \label{xiAdk}
\end{equation}
for $x$ near $o$. Replacing $\phi_i$ by $\phi\phi_i$ for a $K$-invariant $\phi\in C_c^\infty(X)$
with $\phi=1$ near $o$ and $\phi=0$ outside a small neighborhood of $o$,
we may assume (\ref{xiAdk}) holds for all $x\in X$. Then $\phi_i$'s will be called coordinate
functions on $X$ under the basis $\{\xi_i\}$ of $\myp$.

  Similar to the definitions given on $G$, the mean of a random variable $x$ in $X$
or its distribution $\mu$ is defined as $b=e^{\sum_{i=1}^n\mu(\phi_i)\xi_i}o$.
By (\ref{xiAdk}), if $\mu$ is $K$-invariant, then so is $b$.
We will call $x$ or $\mu$ small if $b$ has coordinates $\mu(\phi_i)$,
that is, $\phi_i(b)= \mu(\phi_i)$ for $1\leq i\leq n$.

Any $\xi\in\myg$ is a left invariant vector field on $G$. If $\xi$ is $\Ad(K)$-invariant,
that is, if $\Ad(k)\xi=\xi$ for $k\in K$, then $e^{t\xi}o$ is $K$-invariant and $\xi$ may also
be regarded as a vector field on $X$ given by $\xi f(x)=\frac{d}{dt}f(xe^{t\xi}o)\mid_{t=0}$
for $f\in C_c^\infty(X)$ and $x\in X$, which is $G$-invariant in the sense that $\xi(f\circ g)=
(\xi f)\circ g$ for $g\in G$. In fact, any $G$-invariant vector field on $X$ is given by
a unique $\Ad(K)$-invariant $\xi\in\myp$. Note that if $\xi\in\myg$ is $\Ad(K)$-invariant
and $b\in X$ is $K$-invariant, then $\Ad(b)\xi=\Ad(S(b))\xi$ is $\Ad(K)$-invariant and is independent
of the choice of the section map $S$. By (\ref{xiAdk}), for any $K$-invariant measure $\mu$
on $X$, $\int\mu(dx)\sum_i\phi_i(x)\xi_i$ is $\Ad(K)$-invariant, and so is $\int\mu(dx)\sum_i
\phi_i(x)\Ad(b)\xi_i$ for a $K$-invariant $b\in X$.

Let $\xi,\eta\in\myg$. With a choice of section map $S$, $\xi\eta$ may
be regarded as a second order differential operator on $X$ defined
by $\xi\eta f(x)=\frac{\partial^2}{\partial t\,\partial s}f(S(x)e^{t\xi}e^{s\eta}o)\mid_{t=s=0}$.
It is shown in \cite[part~2]{feinsilver} that $\xi_i\xi_j f(x) =
\frac{\partial^2}{\partial t\,\partial s}f(S(x)e^{t\xi_i+s\xi_j}o)
\mid_{t=s=0}+\sum_{k=1}^n\rho_k^{ij}\xi_kf(x)$ with $\rho_k^{ij}=
-\rho_k^{ji}$. Thus, if $a_{ij}$ is an $n\times n$ symmetric matrix,
then
\begin{equation}
  \sum_{i,j=1}^n a_{ij}\xi_i\xi_jf(x) =
\sum_{i,j=1}^n a_{ij}\frac{\partial^2}{\partial t_i\,\partial t_j}
f(S(x)e^{\sum_p t_p\xi_p}o)\mid_{t_1=\cdots=t_n=0}.
\label{xixif}
\end{equation}
The matrix $a_{ij}$ is called $\Ad(K)$-invariant if $a_{ij}=\sum_{p,q}
a_{pq}[\Ad(k)]_{ip}[\Ad(k)]_{jq}$ for $k\in K$,
where $\{[\Ad(k)]_{ij}\}$ is the matrix representing $\Ad(k)$,
that is, $\Ad(k)\xi_j=\sum_i[\Ad(k)]_{ij}\xi_i$. Then the operator $\sum_{i,j}a_{ij}
\xi_i\xi_j$ is independent of section map $S$ and
is $G$-invariant. In fact, any second order $G$-invariant differential operator $T$
on $X$ with $T1=0$ is such an operator plus a $G$-invariant vector field.
Note that $\sum_{i,j}a_{ij}[\Ad(b)\xi_i]
[\Ad(b)\xi_j]=\sum_{i,j}a_{ij}[\Ad(S(b))\xi_i][\Ad(S(b))\xi_j]$ is
a $G$-invariant operator on $X$ for any $b\in X$ (independent of $S$)
if $a_{ij}$ is $\Ad(K)$-invariant.

A drift or an extended drift $b$, a covariance matrix function $A$,
and a \levy measure function or an extended \levy measure function $\eta$
on $X=G/K$ are defined just as on $G$, in terms of
the coordinate functions $\phi_i$ on $X$, with the additional requirement that
for each $t$, $b_t$ and $\eta(t,\cdot)$ are $K$-invariant, and $A(t)$
is $\Ad(K)$-invariant. With these modifications, \levy triples
and extended \levy triples on $X$ are defined exactly as on $G$. In particular,
for an extended \levy triple $(b,A,\eta)$, $b_{t-}^{-1}b_t=h_t$, where $h_t=
\exp[\sum_{i=1}^n\nu_t(\phi_i)\xi_i]o$ is the mean of the $K$-invariant
probability measure $\nu_t=\eta(\{t\}\times\cdot)+
[1-\eta(\{t\}\times X)]\delta_o$.

A rcll process $x_t$ in $X=G/K$ is said to be represented by an extended \levy triple $(b,A,\eta)$ on $X$ if $x_t=z_tb_t$ and for any $f\in C_c^\infty(X)$,
\begin{eqnarray}
M_tf &=& f(z_t) - \frac{1}{2}
\int_0^t\sum_{j,k=1}^n[\Ad(b_s)\xi_j][\Ad(b_s)\xi_k]
f(z_s)dA_{jk}(s) \nonumber \\
&&\ \ -
\int_0^t\int_X\{f(z_sb_sxb_s^{-1})-f(z_s)-\sum_{j=1}^n\phi_j(x)
[\Ad(b_s)\xi_j]f(z_s)\}\eta^c(ds,dx) \nonumber \\
&&\ \ - \sum_{u\leq
t}\int_X[f(z_{u-}b_{u-}xh_u^{-1}b_{u-}^{-1})-f(z_{u-})] \nu_u(dx)
\label{MtfX}
\end{eqnarray}
is a martingale under the natural filtration $\myF_t^x$ of $x_t$.

By the $K$-invariance of $b_t$ and $\eta(t,\cdot)$, and
the $\Ad(K)$-invariance of $A(t)$, the expression
in (\ref{MtfX}) makes sense. Moreover,
by Taylor expansions of $f(z_sb_sxb_s^{-1})=f(z_se^{\sum_i
\phi_i(x)\Ad(b_s)\xi_i})$ at $x=o$ and $f(z_{u-}b_{u-}xh_u^{-1}b_{u-}^{-1})$
at $x=h_u$, and the properties of $\eta(t,\cdot)$ as an extended \levy measure function,
it can be shown that $M_tf$ given in (\ref{MtfX}) is a bounded random variable.

\begin{theo} \label{thmartrepX}
Let $x_t$ be an inhomogeneous \levy process in $X=G/K$ with $x_0=o$.
Then there is a unique extended \levy triple $(b,A,\eta)$ on $X$ such
that $x_t$ is represented by $(b,A,\eta)$ as defined above.
Moreover, $\eta(t,\cdot)$ is the jump intensity measure of process $x_t$ given
by (\ref{etatBX}). Consequently, $x_t$ is stochastically continuous
if and only if $(b,A,\eta)$ is a \levy triple.

  Conversely, given an extended \levy triple $(b,A,\eta)$ on $X$, there is an inhomogeneous \levy
process $x_t$ in $G$ with $x_0=o$, unique in distribution, such that $x_t$ is
represented by $(b,A,\eta)$.
\end{theo}

\noindent {\bf Proof} \ The theorem may be proved by essentially
repeating the proof of Theorem~\ref{thmartrep} for the corresponding
results on $G$, interpreting a product $xy$ and an inverse $x^{-1}$ on $X=G/K$ as $S(x)y$
and $S(x)^{-1}$ by choosing a section map $S$ as in the preceding discussion, and taking various
functions and sets to be $K$-invariant. As noted in \cite{Liao2009}, the results in \cite{feinsilver}
for stochastically continuous inhomogeneous \levy processes in $G$ hold also on $X$.
In particular, Lemma~\ref{lefein}, which summarizes part of \cite{feinsilver},
holds also on $X$. Note that because $\mu_{ni}=\mu_{t_{n\,i-1},t_{ni}}$ and $\nu_t(\cdot)$
are $K$-invariant, their means $b_{ni}$ and $h_t$ are $K$-invariant, and so is $b_t^n$.
It then follows from (\ref{xiAdk}) that $A^n(t,f)$ is $\Ad(K)$-invariant for $K$-invariant $f\in C_b(G)$,
and (\ref{AtUAtintU}) holds on $X$ with $A(t)$ being a covariance matrix function on $X$. \ $\Box$
\vspace{2ex}

\noindent {\bf Remark 2} \ As in Remark~1 for the representation on $G$, $\eta(t,\cdot)$
in Theorem~\ref{thmartrepX} is independent of the choice for the basis $\{\xi_j\}$ of $\myp$
and coordinate functions $\phi_j$ on $X=G/K$, $A(t)$ is independent of $\{\phi_j\}$
and the $G$-invariant operator $\sum_{j,k=1}^nA(t)\xi_j\xi_k$ on $X$ is independent of $\{\xi_j\}$.
However, $\eta$ may depend on the choice for the origin $o$ in $X$.
\vspace{2ex}

\noindent {\bf Remark 3} \ We note that the representation of inhomogeneous \levy processes
in $G$ by a triple $(b,A,\eta)$ with $b_t$ of finite variation, as stated
in Theorem~\ref{thmartfinvar}, holds also on $X=G/K$ with essentially the same proof,
where $b_i(t)$ are components of $b_t$ under the basis $\{\xi_i\}$ of $\myp$.
\vspace{2ex}

Because the action of $K$ on $X$ fixes $o=eK$, it induces a linear action
on the tangent space $T_oX$ at $o$. The homogeneous space $X=G/K$ is called
irreducible if the $K$-action on $T_oM$ is irreducible, that is,
it has no nontrivial invariant subspace. By the identification of $\myp$
and $T_oX$ via the differential of $\pi$: $G\to X$, $X$ is irreducible if and only if $\myp$ has
no nontrivial $\Ad(K)$-invariant subspace. Because $\pi\circ\exp$: $\myg\to X$ is diffeomorphic
from a neighborhood $V$ of $0$ in $\myp$ onto a neighborhood $U$
of $o$ in $X$, the only $K$-invariant point in $U$ is $o$. For
example, the $n$-dimensional sphere $S^n=O(n+1)/O(n)$ is irreducible.

Let $V$ above be convex. The coordinate functions $\phi_i$
on $X$ may be chosen so that $\sum_{i=1}^n\phi_i(x)\xi_i\in V$
for all $x\in X$. Then the mean $b=
e^{\sum_i\mu(\phi_i)\xi_i}$ of any distribution $\mu$ on $X$ belongs to $U$.
Let $X$ be irreducible. If $\mu$ is $K$-invariant, then so is $b$ and hence $b=o$.
Because the extended drift $b_t$ of an inhomogeneous \levy process in $X$, represented
by $(b,A,\eta)$, is the limit of $b_t^n$ in (\ref{btn}).
it follows that $b_t=o$ for all $t\geq 0$.

As mentioned in \cite{Liao2009}, on an irreducible $X$, up to a constant
multiple, there is a unique $K$-invariant inner product on $T_oX$. This
implies that there is a unique $\Ad(K)$-invariant inner product on $\myp$,
and hence an $\Ad(K)$-invariant matrix is a multiple of the identity matrix.
To summarize, we record the following conclusion.

\begin{prop} \label{prhirredX}
On an irreducible $X=G/K$, any covariance matrix function $A(t)$ is given
by $A(t)=a(t)I$ for some nondecreasing continuous function $a(t)$
with $a(0)=0$, where $I$ is the identity matrix. Moreover, if
the coordinate functions $\phi_i$ are chosen as above, then
the extended drift $b_t$ of an inhomogeneous \levy process in $X$ is trivial,
that is, $b_t=o$ for $t\geq 0$.
\end{prop}

By (\ref{xiAdk}), for $K$-invariant $H\subset X$ and $K$-invariant measure $\mu$
on $X$, $\int_H\mu(dx)\sum_{j=1}^n\phi_j(x)\xi_j$ is a $K$-invariant vector in $\myp$
and hence is $0$ by the irreducibility of $X=G/K$. By Proposition~\ref{prhirredX}, $b_t=o$
and $A_t=a(t)I$, and hence $z_t=x_t$. The integral $\int_0^t\int_X(\cdots)\eta^c(ds,dx)$
in (\ref{MtfX}) is the limit of $\int_0^t\int_{U^c}(\cdots)\eta^c(ds,dx)$ as a $K$-invariant
neighborhood $U$ of $o$ decreases to $o$, and then the third term of the integrand may be
removed. It follows that this integral combined with the sum $\sum_{u\leq t}$ in (\ref{MtfX}) may
be written as $\int_0^t\int_X[f(x_s\tau)-f(x_s)]\eta(ds,d\tau)$, where the integral is
understood as the principal value, that is, as the limit of $\int_0^t\int_{U^c}[\cdots]
\eta(ds,\tau)$ as a $K$-invariant neighborhood $U$ of $o$ shrinks to $o$. By Theorem~\ref{thmartrepX},
we obtain the following simple form of martingale representation on an irreducible $X=G/K$.

\begin{theo} \label{thmartrepirredX}
For an inhomogeneous \levy process $x_t$ in an irreducible $X=G/K$ with $x_0=o$,
there is a unique pair $(a,\eta)$ of a continuous nondecreasing function $a(t)$
with $a(0)=0$ and an extended \levy measure function $\eta(t,\cdot)$ on $X$ such that
\begin{equation}
\forall f\in C_c^\infty(X),\ \ \ \ f(x_t) - \int_0^t\frac{1}{2}\sum_{j=1}^n
\xi_j\xi_jf(x_s)da(s) - \int_0^t\int_X[f(x_s\tau)-f(x_s)]\eta(ds,d\tau) \label{martirredX}
\end{equation}
is a martingale under $\myF_t^x$, where $\int_0^t\int_X[\cdots]\eta(ds,d\tau)$ is
the principal value as described above and is bounded. Moreover, $\eta$ is
the jump intensity measure of $x_t$ given by (\ref{etatBX}). Consequently, $x_t$ is stochastically continuous
if and only if $\eta$ is a \levy measure function.

  Conversely, given a pair $(a,\eta)$ as above, there is an inhomogeneous \levy process $x_t$
in $X$ with $x_0=o$, unique in distribution, such that (\ref{martirredX}) is a martingale under $\myF_t^x$.
\end{theo}

An inhomogeneous \levy process $g_t$ in $G$ is called $K$-conjugate invariant
if its transition function $P_{s,t}$ is $K$-conjugate invariant, that is,
if $P_{s,t}(f\circ c_k)=(P_{s,t}f)\circ c_k$ for $f\in\myB_+(G)$ and $k\in K$,
where $c_g$: $x\mapsto gxg^{-1}$ is the conjugation map on $G$. Note that
when $g_0=e$, this is equivalent to saying that the process $g_t$ has
the same distribution as $kg_tk^{-1}$ for any $k\in K$.

\begin{theo} \label{thxtgto}
Let $g_t$ be a $K$-conjugate invariant inhomogeneous \levy process in $G$
with $g_0=e$. Then $x_t=g_to$ is an inhomogeneous \levy process in $X=G/K$. Conversely,
if $x_t$ is an inhomogeneous \levy process in $X$ with $x_0=o$, then
there is a $K$-conjugate invariant inhomogeneous \levy process $g_t$ in $G$ such
that processes $x_t$ and $g_to$ are equal in distribution.
\end{theo}

\noindent {\bf Proof} \ Let $g_t$ be a $K$-conjugate invariant inhomogeneous \levy
process in $G$ with $g_0=e$ and let $x_t=g_to$. Then for $s<t$ and a section map $S$
on $X$, $S(x_s)=g_sk_s$ for some $\myF_s^g$-measurable $k_s\in K$, and $S(x_s)^{-1}x_t=
k_s^{-1}g_s^{-1}g_to\stackrel{d}{=}g_s^{-1}g_to$ (equal in distribution because
of $K$-conjugate invariance of $g_t$). This shows that $S(x_s)^{-1}x_t$ is independent
of $\myF_s^g$ and its distribution is independent of $S$. By Proposition~\ref{prii}, $x_t$
is an inhomogeneous \levy process in $X$.

Now let $x_t$ be an inhomogeneous \levy process in $X$ with $x_0=o$, represented by
an extended \levy triple $(b,A,\eta)$ on $X$. We will define an extended \levy
triple $(\bar{b},\bar{A},\bar{\eta})$ on $G$ from $(b,A,\eta)$. Choose a section map $S$
on $X$ such that $S(x)=e^{\sum_{i=1}^n\phi_i(x)\xi_i}o$ for $x$ contained in
a $K$-invariant neighborhood $U$ of $o$. By (\ref{xiAdk}),
\begin{equation}
\forall k\in K\ {\rm and}\ x\in U,\ \ \ \ kS(x)k^{-1} = S(kx). \label{kSxkinv}
\end{equation}

The basis $\{\xi_1,\ldots,\xi_n\}$ of $\myp$ is now extended to be
a basis $\{\xi_1,\ldots,\xi_d\}$ of $\myg$ so that $\{\xi_{n+1},\ldots,
\xi_d\}$ is a basis of $\myk$. Let $\psi_1,\ldots,\psi_d\in C_c^\infty(G)$ be
the associated coordinate functions on $G$, that is, $g=e^{\sum_{i=1}^d\psi_i(g)\xi_i}$
for $g$ near $e$. Because $kgk^{-1}=e^{\sum\psi_i(g)\Ad(k)\xi_i}$,
\begin{equation}
\forall k\in K,\ \ \ \ \sum_{i=1}^d(\psi_i\circ c_k)\xi_i = \sum_{i=1}^d\psi_i[\Ad(k)\xi_i]
\label{psiAdkxi}
\end{equation}
holds near $e$. Replacing $\psi_i$ by $\psi\psi_i$ for a $K$-conjugate invariant $\psi\in C_c^\infty(G)$ such
that $\psi=1$ near $e$ and $\psi=0$ outside a small neighborhood of $e$,
we may assume (\ref{psiAdkxi}) holds on $G$.

Because for $x$ near $o$, $e^{\sum_{i=1}^n\phi_i(x)\xi_i}=S(x)=
e^{\sum_{i=1}^d\psi_i(S(x))\xi_i}$, it follows that
\begin{equation}
\psi_i\circ S = \phi_i,\ \ 1\leq i\leq n,\ \ \ \ {\rm and} \ \ \ \
\psi_i\circ S = 0,\ \ n+1\leq i\leq d \label{psiSphin}
\end{equation}
holds near $o$. The $\psi_i$ may be modified outside a neighborhood of $e$ so
that (\ref{psiSphin}) holds on $X$ and (\ref{psiAdkxi}) still holds on $G$. This may be accomplished
by replacing $\psi_i$ by $(\phi\circ\pi)\psi_i+(1-
\phi\circ\pi)(\phi_i\circ\pi)$ for $i\leq n$ and by $(\phi\circ\pi)
\psi_i$ for $i>n$, where $\phi$ is a $K$-invariant smooth function
on $X$ such that $\phi=1$ near $o$ and $\phi=0$ outside a small
neighborhood of $o$.

Let $\bar{h}_t=e^{\sum_j\nu_t(\phi_j)\xi_j}$. Then $\pi(\bar{h}_t)=b_{t-}^{-1}b_t$.
By (\ref{xiAdk}), $\bar{h}_t$ is $K$-conjugate invariant.
There is a partition $0=t_0<t_1<t_2<\cdots<t_n\uparrow\infty$ such that $b_{t_i}^{-1}b_t
\in U$ for $t\in [t_i,\,t_{i+1})$. Define $\bar{b}_t= S(b_t)$
for $t<t_1$ and $\bar{b}_{t_1}=\bar{b}_{t_1-}\bar{h}_{t_1}$, and inductively let $\bar{b}_t=\bar{b}_{t_i}S(b_{t_i}^{-1}b_t)$
for $t_i\leq t<t_{i+1}$ and $\bar{b}_{t_{i+1}}=\bar{b}_{t_{i+1}-}\bar{h}_{t_{i+1}}$.
Then $\bar{b}_t$ is a $K$-conjugate invariant extended drift in $G$ with $\pi(\bar{b}_t)=b_t$.
Let $\bar{A}(t)$ be the $d\times d$ matrix given by  $\bar{A}_{ij}(t)=A_{ij}(t)$ for $i,j\leq n$
and $\bar{A}_{ij}=0$ otherwise. Then $\bar{A}(t)$ is a covariance matrix function on $G$
that is $\Ad(K)$-invariant, that is, $[\Ad(k)]\bar{A}(t)[\Ad(k)]'=\bar{A}(t)$ for $k\in K$,
where $[\Ad(k)]$ is the matrix representing the linear map $\Ad(k)$: $\myg\to\myg$
under the basis $\{\xi_1,\ldots,\xi_d\}$. For $f\in\myB_+(G)$, let
\begin{equation}
\bar{\eta}(t,f) = \int_X\int_K f(kS(x)k^{-1})dk\eta(t,dx).
\label{baretat}
\end{equation}
Then $\bar{\eta} (t,\cdot)$ is a $K$-conjugate invariant measure function. Its continuous
part is $\bar{\eta}^c(t,f)=\int\int f(kS(x)k^{-1})dk\eta^c(t,dx)$ and for each $t>0$, $\bar{\nu}_t=
\bar{\eta}(\{t\},\cdot)+[1-\bar{\eta}(\{t\},G)]\delta_e$ is given by $\bar{\nu}_t(f)=
\int\int f(kS(x)k^{-1})dk\nu_t(dx)$.

By (\ref{xiAdk}) and (\ref{psiSphin}) on $X$, and (\ref{psiAdkxi}) on $G$, $\sum_{j=1}^d\bar{\nu}_t
(\psi_j)\xi_j=\sum_{j=1}^d\int_X\int_K[\psi_j(kS(x)k^{-1})\xi_j]dk\nu_t(dx)=\sum_{j=1}^d
\int_X\int_K\{\psi_j(S(x))[\Ad(k)]\xi_j\}dk\nu_t(dx)=\sum_{j=1}^n\int_X\{\int_K\phi_j(x)[\Ad(k)]\xi_j\}dk
\nu_t(dx)=\sum_{j=1}^n\nu_t(\phi_j)\xi_j$. This shows that $\bar{h}_t=e^{\sum_j\nu_t(\phi_j)\xi_j}$ defined
earlier is the mean of $\bar{\nu}_t$.

We now show that $\bar{\eta}$ is an extended \levy measure function on $G$
by checking $\bar{\eta}^c(t,\Vert\psi_\cdot\Vert^2)<\infty$ and $\sum_{u\leq t}
\bar{\nu}_u(\Vert\psi_\cdot-\psi_\cdot(\bar{h}_u)\Vert^2)<\infty$, where $\Vert\psi_\cdot\Vert$ is
the Euclidean norm of $\psi_\cdot=(\psi_1,\ldots,\psi_d)$.
Let $V=\pi^{-1}(U)$. Then $\bar{\eta}(t,V^c)=
\eta(t,U^c)<\infty$. By (\ref{kSxkinv}) and (\ref{psiSphin}), the first condition follows
from $\bar{\eta}^c(t,\Vert\psi_\cdot\Vert^21_V)=\eta^c(t,\Vert\phi_\cdot\Vert^21_U)<\infty$.
To verify the second condition, we may assume all $\bar{\nu}_t$ are small,
then $\psi_i(\bar{h}_t)=\bar{\nu}_t(\psi_i)$ and the second condition may be
written as $\sum_{u\leq t}\int_G\bar{\nu}_u(dx)\Vert\int_G\bar{\nu}_u(dy)
[\psi_\cdot(x)-\psi_\cdot(y)]\Vert^2<\infty$. Because $\sum_{u\leq t}\bar{\nu}_u(V^c)=\sum_{u\leq t}
\nu_u(U^c)<\infty$, this is equivalent to $\sum_{u\leq t}\int_V\bar{\nu}_u(dx)\Vert\int_V\bar{\nu}_u(dy)
[\psi_\cdot(x)-\psi_\cdot(y)]\Vert^2<\infty$, but this is the same as $\sum_{u\leq t}\int_U\nu_u(dx)
\Vert\int_U\nu_u(dy)[\phi_\cdot(x)-\phi_\cdot(y)]\Vert^2<\infty$, which holds by (\ref{sumnu}).

Now let $g_t$ be an inhomogeneous \levy process in $G$ with $g_0=e$, represented
by the extended \levy triple $(\bar{b},\bar{A},\bar{\eta})$. Write $g_t=
\bar{z}_t\bar{b}_t$ and for $f\in C_c^\infty(G)$, let $\bar{M}_tf(\bar{z}_\cdot)$
be the martingale given in (\ref{Mtf}) with $(b,A,\eta)$ and $\phi_i$ replaced by $(\bar{b},
\bar{A},\bar{\eta})$ and $\psi_i$. By the $K$-conjugate invariance of $\bar{b}$ and $\bar{\eta}$,
the $\Ad(K)$-invariance of $\bar{A}$, and (\ref{psiAdkxi}), it is easy to show that
for $k\in K$, $\bar{M}_t(f\circ c_k)(\bar{z}_\cdot)=\bar{M}_tf(k\bar{z}_\cdot k^{-1})$.
This means that $kg_tk^{-1}=k\bar{z}_tk^{-1}\bar{b}_t$ is also represented by $(\bar{b},
\bar{A},\bar{\eta})$. By the unique distribution in Theorem~\ref{thmartrep}, $g_t$
and $kg_tk^{-1}$ have the same distribution as processes, and hence $g_t$ is $K$-conjugate invariant.

Let $x_t'=g_to$. By the first part of Theorem~\ref{thxtgto} which has been proved, $x_t'$ is an inhomogeneous \levy process in $X$.
We may write $x_t'=z_tb_t$ with $z_t=\bar{z}_to$. By the construction of $(\bar{b},\bar{A},\bar{\eta})$ from $(b,A,\eta)$, it can be shown
that $\bar{M}_t(f\circ\pi)(\bar{z}_\cdot)$ is the martingale in (\ref{MtfX}). This shows that $x_t'$ is represented by $(b,A,\eta)$. By the unique distribution
in Theorem~\ref{thmartrepX}, $x_t'$ and $x_t$ are equal in distribution. \ $\Box$
\vspace{2ex}

\noindent {\bf Acknowledgment} \ Helpful comments from David Applebaum and an anonymous referee have led to a considerable improvement of the exposition.

\noindent Department of Mathematics, Auburn University,
Auburn, AL 36849, USA.

\noindent Email: liaomin@auburn.edu

\end{document}